\documentclass[12pt]{article}
\usepackage{amsmath,amssymb}

\def\CC{{\mathbb C}}
\def\ZZ{{\mathbb Z}}

\def\RR{{\mathbb R}}
\def\QQ{{\mathbb Q}}
\def\DD{{\mathbb D}}
\def\PP{{\mathbb P}}

\def\MM{{\mathbb M}}
\def\ZZ{{\mathbb Z}}

\newcommand {\eps} {\epsilon}
\newcommand {\BOX} {\hfill \rule{2mm}{2mm} \bigskip}

\newcommand {\be} {\begin{equation}}
\newcommand {\eeqn} {\end{equation}}
\newcommand {\bea} { \begin{eqnarray}}
\newcommand {\eea} {\end{eqnarray}}
\newcommand {\beas} { \begin{eqnarray*}}
\newcommand {\eeas} {\end{eqnarray*}}
\newcommand {\ra} {\rightarrow}

\newcommand {\ol} {\overline}

\newtheorem {lemma} {LEMMA} [section]
\newtheorem {theorem}[lemma]{THEOREM} 
\newtheorem {prop}[lemma]{PROPOSITION}
\newtheorem {cor}[lemma]{COROLLARY}
\newtheorem {definition}[lemma]{DEFINITION}
\newtheorem {question}[lemma]{QUESTION}

\numberwithin{equation}{section}

\begin{document}

\title{Algebraic Surfaces Holomorphically Dominable by $\CC^2$}
\author{Gregery T. Buzzard\thanks{Partially supported by an NSF grant.} $\ $
and Steven S. Y. Lu\thanks{Partially supported by MSRI and NSERC.}}
\date{}
\maketitle


\section{Introduction}

An $n$-dimensional complex manifold $M$ is said to be (holomorphically) 
{\it dominable} by $\CC^n$ if
there is a map $F:\CC^n \ra M$ which is holomorphic such that the Jacobian
determinant $\det(DF)$ is not identically zero.  Such a map $F$ is
called a {\it dominating map}.  In this paper, we attempt to classify 
algebraic surfaces $X$ which are dominable by $\CC^2$ using a combination of
techniques from algebraic topology, complex geometry and analysis. One
of the key tools in the study of algebraic surfaces is the notion of
Kodaira dimension (defined in section 2).  By 
Kodaira's pioneering work \cite{K} and its extensions 
(see, for example, \cite{CG} and \cite{KO}), an algebraic surface
which is dominable by $\CC^2$ must have Kodaira dimension less than
two. Using the Kodaira dimension and the fundamental group of $X$, we
succeed in classifying algebraic surfaces which are dominable by
$\CC^2$ except for certain cases in which $X$ is an algebraic surface 
of Kodaira dimension zero and the case when $X$ is rational without
any logarithmic $1$-form. More specifically, in the case when $X$ is 
compact (namely projective), we need to exclude only the case when $X$ is
birationally equivalent to a K3 surface (a simply connected compact complex
surface which admits a globally non-vanishing holomorphic 2-form)
that is neither elliptic nor Kummer (see sections 3 and 4 for the
definition of these types of surfaces).  

With the exceptions noted above, we show that for any algebraic
surface of Kodaira dimension less than 2, dominability by $\CC^2$ is 
equivalent to the apparently weaker requirement of the existence of a
holomorphic image of $\CC$ which is Zariski dense in the surface. With
the same exceptions, we will also show the very interesting and 
revealing fact that dominability by $\CC^2$ is preserved even if a
sufficiently small neighborhood of any finite set of points is removed
from the surface. In fact, we will provide a complete classification in the 
more general category of (not necessarily algebraic) compact complex surfaces 
before tackling the problem in the case of non-compact algebraic surfaces.

We remark that both elliptic K3 and Kummer K3 surfaces are dense in the
moduli space of K3 surfaces; the former is dense of codimension-one
while the latter is dense of codimension sixteen in this moduli space
(see \cite{PS,LP}) and intersects the former transversally (these
density results hold also in any universal family). 
Dominability by $\CC^2$ holds for both types of K3 surfaces. This
suggests that it might hold for all K3 surfaces so that our statements
above would be valid without exception for projective (and, more generally, 
for compact K\"ahler) surfaces.  Indeed, their density plus Brody's
lemma (\cite{Br}) tells us that every K3 surface contains a
non-trivial holomorphic image of $\CC$ and that the generic 
K3 surface, which is non-projective but remains K\"ahler,
even contains such an image that is Zariski dense. We mention here that
dominability by $\CC^2$ can be shown for some non-elliptic K3 surfaces
which are close to Kummer surfaces
using an argument similar to that of section 6; for length
considerations, we omit this non-elliptic case from this paper.  
However, we note that the statement equating dominability 
to the weaker condition of having a Zariski dense image of $\CC$
is quite false in the non-K\"ahler category, as is amply demonstrated by
Inoue surfaces (see \cite{In0} or \cite[V.19]{BPV}).

Observe that if there is a dominating map $F:\CC^2 \ra X$, then there is
also a holomorphic image of $\CC$ which is Zariski dense:
First we may assume that the Jacobian of $F$ is non-zero at
the origin.  Defining $h:\CC \ra \CC^2$ by $h(z) = (\sin(2 \pi z),
\sin(2\pi z^2))$, we see that $h(n)=(0,0)$ with corresponding tangent
direction $(2\pi, 4\pi n)$ for each $n \in \ZZ$.  Taking $F \circ h$, we
obtain a holomorphic image of $\CC$ with an infinite number of tangent
directions at one point, which implies that the image is Zariski dense.

We say that an algebraic variety $X$ satisfies property C
if every holomorphic image of $\CC$ in $X$ is algebraically
degenerate; i.e., is not Zariski dense.  Our first main result is that
for algebraic surfaces of Kodaira dimension less than 2 and with the
exceptions mentioned above, dominability by $\CC^2$
is equivalent to the failure of property C. We will
state only the main results in the projective category in this introduction
for simplicity but will discuss fully the compact non-projective case and 
much of the quasi-projective case in this paper.
\begin{theorem} \label{thm:1}
Let $X$ be a projective surface of Kodaira dimension less than $2$ and
suppose that $X$ is not birational to a K3 surface which is either
elliptic or Kummer.  Then $X$ is dominable
by $\CC^2$ if and only if it does not satisfy property C.
Equivalently, there is a  
dominating holomorphic map $F:\CC^2 \ra X$
if and only if there is a holomorphic image of
$\CC$ in $X$ which is Zariski dense.
\end{theorem}

By a recent result of the second named author, this theorem is also
true for a projective surface of Kodaira dimension 2, which is the
maximum Kodaira dimension for surfaces.  As previously mentioned, a
surface of Kodaira dimension 2 is not dominable by $\CC^2$ \cite{K};
indeed, a surface of Kodaira dimension 2 is precisely a surface which
admits a possibly degenerate hyperbolic volume form. Thus in the case
of Kodaira dimension 2, theorem~\ref{thm:1} can be established by
showing that such a surface satisfies property C.  The question of
whether a variety of maximum Kodaira dimension satisfies property C
was first raised explicitly by Serge Lang \cite{Lang}.

\smallskip  

In the following theorem we give, again modulo the above mentioned
exceptions, a classification of 
projective surfaces which are dominable by $\CC^2$ and hence a
classification of projective surfaces of Kodaira dimension less than 2
which fail to satisfy property C.  We will do this in terms
of the Kodaira dimension and the fundamental group, both of which are
invariant under birational maps.

\begin{theorem}  \label{thm:2}
A projective surface $X$ not birationally equivalent to a K3 surface
is dominable by $\CC^2$ if and only if 
it has Kodaira dimension less than two and its fundamental group
is a finite extension of an abelian group (of even rank four or less).
If $\kappa (X)=-\infty$, then the fundamental group condition can be
replaced by the simpler condition of non-existence of more than one 
linearly independent holomorphic one-form. If $\kappa (X)=0$ and $X$
is not birationally equivalent to a K3 surface, then $X$ is dominable
by $\CC^2$. If $X$ is birationally equivalent to an elliptic K3
surface or to a Kummer K3 surface, then $X$ is dominable by $\CC^2$.
\end{theorem}

As with theorem \ref{thm:1}, this theorem fails if we include compact
non-K\"ahler surfaces (even after simple minded modification of this
theorem). For instance, the Kodaira surfaces are dominable by $\CC^2$ but
their fundamental groups are not finite extensions of abelian groups
(\cite{Ko4}). But this theorem remains valid in the K\"ahler category,
thanks, for example, to Kodaira's result that all K\"ahler surfaces
are deformations of projective surfaces (\cite{Ko2}, \cite{Ko3}). 

More general versions of theorem~\ref{thm:1} and theorem~\ref{thm:2}
for compact complex surfaces will be given at the end of section 4.

\smallskip

In the quasi-projective category, we also prove the analogue of
theorem~\ref{thm:1} modulo the same exceptions mentioned in the
beginning, following mainly the work of Kawamata \cite{K1} and
M. Miyanishi \cite{M}.  In this setting, the analogue of the fundamental
group characterization requires the study of a new but very
natural class of objects of complex dimensional one that are related to 
orbifolds. As for explicit examples, we will work out 
theorems~\ref{thm:1} and an analogue of~\ref{thm:2} for
the complement of a reduced curve $C$ in $\PP^2$ in the case when
$C$ is normal crossing, where we show that
dominability is characterized by $\deg C\leq 3$, and for the
overlapping case in which $C$ is either not a rational curve of high
degree or has at most one singular point. Here, the most fascinating
and revealing example is the case in which $C$ is a non-singular cubic
curve, whose complement is a noncompact analogue of a K3 surface. The
question of the dominability of the complement of a 
non-singular cubic was discussed by Bernard Shiffman at MSRI in 1996,
and the positive resolution of this problem served as the first result
in and inspiration for this paper.  

\smallskip

The key tools we introduce
here for constructing dominating maps are the mapping theorems
we establish via a combination of complex geometry and analysis. One of these
theorems utilizes Kodaira's theory of Jacobian fibrations to deal with
general elliptic fibrations (see section 3). Other such theorems construct
the required self-maps of $\CC^2$ directly via complex analysis
to deal with $\CC^\ast$-fibrations, abelian and Kummer surfaces.

\smallskip
In particular, the constructions in sections 4 and 6 show that given any
complex 2-torus and any finite set of points in this torus, there is
an open set containing this finite set and a dominating map from
$\CC^2$ into the complement of the open set.  This should be compared
with \cite{green} in which it was claimed that the complement of any
open set in a simple complex torus is Kobayashi hyperbolic (a complex
torus is simple if it has no nontrivial complex subtori).  There is no 
contradiction because it was later realized that the proof given in
\cite{green} is incorrect since the topological closure of a
complex one-parameter group need not be a complex torus.   Despite
this, the validity of this claim appears to have been an open question
until the current paper, which shows the claim to be false in
dimension 2.  The $n$-dimensional analogue of our result is given in
\cite{Bu}. 

\smallskip

	Many of the tools and results we develop may be of interest to other 
areas of mathematics besides complex analysis and holomorphic geometry,
especially to Diophantine (arithmetic) geometry in view of the connection
between the transcendental holomorphic properties and arithmetic properties
of algebraic varieties.
For example, the important technique of constructing sections of elliptic
fibrations, which is very difficult to achieve in the algebro-geometric
category but certainly useful in arithmetic  and algebraic geometry,
turns out to be quite natural and relatively easy to do in the
holomorphic category. Also, we undertake a global study, from the
viewpoint of holomorphic geometry, of the monodromy action on the
fundamental group of an elliptic fibration. 
Needless to say, without the deep and beautiful contributions of Kodaira
on complex analytic surfaces, we would not be able to go much beyond
dealing with some special examples, as 
is the case with much of the scarce literature on the subject. However,
we have not avoided, due to the nature of this joint paper,
giving elementary lemmas and proofs while avoiding the unnecessary
full force of Kodaira's theory on elliptic fibrations, especially as 
we deal with fibrations over curves that are not necessarily quasiprojective.

\smallskip

The paper is organized as follows.  Section 2 introduces some basic
birational invariants and general notation and provides a list
of the classification of projective surfaces. 
Section 3 deals with projective surfaces not of zero Kodaira dimension
and solves the dominability problem completely for elliptic
fibrations, including the non-algebraic ones. Section 4 deals with the
remaining projective and compact complex cases while section 5 deal
with the non-compact algebraic surfaces.  Section 6 goes beyond these
theorems to deal with algebraic surfaces minus small open balls. 

\smallskip

We are very grateful to Bernard Shiffman for posing the question which
motivated and inspired this paper and for his constant encouragement
during its preparation.


\section{Classification of algebraic surfaces} %
\label{sec:classification}

In this section we will first introduce some basic invariants in the 
(logarithmic) classification theory of
algebraic varieties (see \cite{Ii} for more details,
also compare with \cite{Ue}). Then we will
provide a list of the birational classification of projective
surfaces and discuss briefly the dominability problem in the
quasi-projective category. Finally, we will introduce the more general
category of compactifiable complex manifolds and a basic invariant
which distinguishes the algebraic case in dimension two.
\medskip

Let $\bar{X}$ be a complex manifold with a normal
crossing divisor $D$. This means that around any point $q$ of $\bar{X}$,
there exist a local coordinate $(z_1,...,z_n)$ 
centered at $q$ such that, for some $r \leq n$, $D$ is defined by
$z_1z_2...z_r=0$ in this coordinate neighborhood.
If all the components of $D$ are  smooth, then $D$ is called a simple
normal crossing divisor.
Following Iitaka (\cite{Ii}), we define the logarithmic
cotangent sheaf $\Omega_{\bar{X}}( {\rm log} D)$ as the locally free 
subsheaf of the sheaf
of meromorphic 1-forms, whose restriction to 
$X=\bar X\setminus D$ is identical to $\Omega_X$
and whose localization at any point $q\in D$ is given by
$$\Omega_{\bar{X}}( {\rm log} D)=
\sum_{i=1}^r {\cal O}_{\bar X,q}\frac{dz_i}{z_i} +
\sum_{j=r+1}^n {\cal O}_{\bar X,q} dz_j,$$
where the local coordinates $z_1,...,z_n$ around $q$ are chosen as before.
Its dual, the logarithmic tangent sheaf $T_{\bar{X}}(- {\rm log} D)$,
is a locally free subsheaf of
$T_{\bar{X}}$. We will follow a general abuse of notation and use the
same notation to denote both a locally free sheaf and a vector bundle. 

By an algebraic variety in this paper,
we mean a complex analytic space $X_0$ such that $X_0$
has an algebraic structure in the following sense:
$X_0$ is covered by a finite number of neighborhoods, each of which is
isomorphic to a closed analytic subspace of a complex vector space defined by
polynomial equations and which piece together with rational coordinate
transformations. A proper birational map from $X_0$ to another variety $X_1$
is, by the graph definition, an algebraic subvariety of $X_0\times X_1$ which
projects generically one-to-one onto  each factor as a proper
morphism. If such a map exists, 
we say that the two varieties are properly birational.
This notion corresponds to that of a bimeromorphic map in the holomorphic
context. Two algebraic varieties are said to be birationally equivalent if
they have isomorphic rational function fields; or equivalently, if they have
birational compactifications.
Hironaka's resolution of singularity theorem
\cite{Hi} (an elementary proof of which can be found in \cite{BM})
implies that given any algebraic variety $X_0$, there is a smooth
projective variety $\bar{X}$ with a simple normal crossing divisor $D$ such
that $X = \bar{X} \setminus D$ is properly birational to $X_0$. If $X_0$
is smooth, then we can even take $X$ to be $X_0$ so that $X_0$ can be 
compactified by adding a simple normal crossing boundary divisor. 
In this paper, a surface will mean
a complex two dimensional manifold while a curve that is not  
explicitly a subvariety (or a subscheme) will mean a (not necessarily
quasi-projective) complex one-dimensional manifold. All surfaces and curves
are assumed to be connected.  In particular, every algebraic surface
is isomorphic to the complement of a  finite  set of  transversely
intersecting smooth curves without triple intersection in some
projective surface. We will use the Enriques-Kodaira classification of
compact surfaces to simplify our problem for surfaces.

\medskip

One of the most important invariants under proper birational maps 
is the (logarithmic) Kodaira dimension.  
Let $X_0$, $X$, $\bar{X}$, and $D$ be as above, and
let $K_{\bar{X}} = \det_\CC (T_{\bar{X}}^\vee)$ where $T_{\bar{X}}^\vee$
is the complex cotangent bundle of $\bar X$. The (holomorphic) line
bundle $K_{\bar X}$ is called the canonical bundle of $\bar{X}$. 
Identifying a line
bundle and its sheaf of holomorphic sections, we define a new line bundle
$K=K_{\bar{X}}(D)=K_{\bar X}\otimes {\cal O}(D)$ corresponding to the sheaf of
meromorphic sections of $K_{\bar X}$ 
which are holomorphic except for simple poles along
$D$ (see Griffiths and Harris \cite{GH} among many other standard references).
In fact, $$K=\det \Omega_{\bar X} (\log D).$$ 
This line bundle on $\bar X$ is called the
logarithmic canonical bundle of $X=\bar X\setminus D$, or more specifically,
of $(\bar X, D)$. We will write tensor products of line bundles additively
by a standard abuse of notation; for example, $mK=K^{\otimes m}$.
Given a projective manifold $\bar Y$ and a birational
morphism $f:\bar Y \ra \bar X$ such that $f^{-1}(D)$ is the same as
a normal crossing divisor $E$ in $\bar Y$, then any section of  $mK$ 
as a tensor power of rational 2-form on $X$ pulls
back via $f$ to a section of $mK_{\bar Y}(E)$. Conversely, any section of
$mK_{\bar Y}(E)$ pulls back (via $f^{-1}$) to a section of $mK$ outside
a codimension-two subset 
(the indeterminacy set of $f^{-1}$ ),
which therefore extends to a section of $mK$ by the
classical extension theorem of Riemann. It follows that,
for every positive integer $m$, 
$h^0(mK):=\dim H^0(mK)$  
is independent of the choice of $\bar X$ for $X_0$ and is a
proper birational invariant of $X_0$. This allows us to introduce the
following birational invariant of $X_0$.

\begin{definition}
The Kodaira dimension of $X_0$ is defined as
$$
\bar{\kappa} (X _{0}) = \limsup_{m \ra \infty} \frac{\log h^0(mK)}{\log m}.
$$
\end{definition} 
The simpler notation $\kappa(X_0)$ is used when $X_0$ is projective.
The Riemann-Roch formula shows that $\bar{\kappa}(X_0)$ 
takes values in the set
$$ \{-\infty,0,1,...,\dim X_0\}.$$
By the same argument as that for $h^0(mK)$, we see that another proper
birational invariant is given by the 
(logarithmic) irregularity of $X_0$ defined by
$$\bar q(X_0)=h^0(\Omega_{\bar X}(\log D)).$$
If $D=0$, then $\bar q(X_0)$ 
is just the dimension of the space of global holomorphic
one-forms $q(X)=h^0(\Omega_X)$ on $X$.

If $\bar{\kappa}(X_0) = \dim(X_0)$, then $X_0$
is called a variety of general type. A theorem of Carlson
and Griffiths \cite{CG} (see also Kodaira \cite{K})
says that $X_0$ cannot be dominated 
(even meromorphically) by $\CC^n$ in this case.
Hence for both theorem~\ref{thm:1} and
theorem~\ref{thm:2}, we need consider only those surfaces with Kodaira
dimension less than 2.

A projective surface $X$ whose canonical bundle has non-negative intersection
with (or, equivalently, non-negative degree when restricted to)
any curve in $X$ is called minimal. We say that $K_X$ is
nef (short for numerically effective) in this case. In general, we say that
a line bundle $L$ on $X$ is nef if $L\cdot C \geq 0$ for any curve $C$ in $X$.

Every algebraic surface is either projective or
admits a projective compactification by adding a set of smooth curves
with at most normal crossing singularities.  Moreover, the
Enriques-Kodaira classification  \cite[Ch. VI]{BPV} 
says that a projective surface admits a
birational morphism (as a composition of blowing up
smooth points) to one of the following.  
\begin{itemize}
\item[\rm (0)] A surface of general type:  $\kappa=2$.
\item[\rm (1)]
$\PP^2$ or a ruled surface over a curve $C$ of genus $g=h^0(\Omega_C)$ 
(that is, a holomorphic
$\PP^1$ bundle over $C$). The latter is birationally equivalent to
$\PP^1\times C$.  Here,  ${\kappa}=-\infty$.
\item[\rm (2)]
An abelian surface (a projective torus given by $\CC ^ {2}$/a lattice).
Here, $\kappa=0$.
\item[\rm (3)]
A K3 surface (a simply connected surface with trivial canonical bundle).
$\kappa=0$.
\item[\rm (4)]
A minimal surface with the structure of an elliptic fibration (see
section~\ref{subsec:elliptic}). \\ Here ${\kappa}$ can be $0$, $1$,
or $-\infty$.
\end{itemize}

The characteristic property of the surfaces listed above is the absence 
of $(-1)$-curves.  A $(-1)$-curve is a smooth
rational curve (image of $\PP^1$) 
in a surface with self-intersection $-1$, i.e. whose
normal bundle has degree $-1$. From Castelnuovo's criterion 
\cite[III4.1]{BPV}, a $(-1)$ curve
is always the blow-up of a (smooth) point on a surface. 
A simple argument (via the linear independence of 
the total transform of blown up $(-1)$-curves in $H_1$)
shows that, given any projective surface with $\kappa<2$, 
one can always reach one of the surfaces listed above by 
blowing down $(-1)$-curves a finite number of times. 
It is a standard
fact that a projective surface with $\kappa\geq 0$ is { minimal} if and only
if it does not have any $(-1)$-curve, and that such a surface is the unique
one in its birational class having this property.
\smallskip

Let $X_0$ be an algebraic surface having a compactification $\bar{X}_0$
which is birational to one of the model surfaces listed above, say $\bar X$.
There is a maximum  Zariski open subset $U$ of  $X_0$ that is properly
birational to the complement of a reduced divisor $C$ and 
a finite set $T'$ of points in $\bar X$.  Now the 
indeterminacy set of this proper birational map from 
$X=\bar X\setminus \{C\cup T\}$
to $U\subseteq X_0$ must
consist of a finite set of points. So to produce a dominating map from
$\CC^2$ to $X_0$, it suffices to produce, for each finite
set of points $T$ in $\bar X$, a dominating map
from $\CC^2$ into the complement of $T$ in $\bar X\setminus C$.
Nevertheless, $\bar X\setminus C$ may not be dominable by $\CC^2$
when $X_0$ is dominable by $\CC^2$; for example, a point on $X_0$ may 
correspond to an infinitely near point on $\bar X$ over a point of $C$. 
However, if we
think of $X_0$ as an open subset of the space of
infinitely near points of $\bar X$, then we can recover
the equivalence in dominability through the above procedure (see section 5).

\medskip

Although we have chosen to introduce and state our results 
so far in the algebraic
category for simplicity, we will in fact 
deal with a more general class of surfaces in the next
two sections: the class of compactifiable surfaces. These are Zariski open
subsets of compact complex surfaces and the invariants $\bar \kappa$
and $\bar q$ carry over to them verbatim as they are defined by 
compactifications with normal crossing divisors, which exist by 
complex surface theory.
If a surface $X$ is compact, the transcendency degree 
$a(X)$ of the field of meromorphic functions on $X$ is, by definition, a 
bimeromorphic invariant and is called the algebraic dimension.


\section{Compact surfaces with $\kappa\neq 0$ and $a\neq 0$}  
\label{sec:compact}

In this section we solve the $\CC^2$ dominability problem
for compact surfaces whose Kodaira dimension and algebraic dimension are
both non-zero. The bulk of this section is devoted to 
the case of elliptic fibrations, which we treat completely, 
including all the noncompact cases.
In particular, we solve our problem  for every projective
surface that is birational to a minimal one listed in (1) and (4) above.
Cases (2) and (3) will be discussed in section~\ref{sec:noncompact}. 


\subsection{Projective surfaces with Kodaira dimension $-\infty$}

Since any $\PP^1$-bundle over a curve $C$ is birational to the trivial
$\PP^1$-bundle over $C$ and since $\PP^2$ is birational to $\PP^1\times\PP^1$,
any projective surface $X$ with $\bar{\kappa}(X)=-\infty$ is birational
to a surface $Y$ which is a trivial $\PP ^ {1}$-bundle over a curve
$C$ of genus $g:=h^0(\Omega_C)$.  In the case where $C$  is of genus
$g>1$, any holomorphic image of $\CC$  in $Y$ must lie in a fiber of
the bundle since  $C$  is hyperbolic.  Hence $X$ satisfies 
property C and so cannot be dominated by $\CC^2$.
In the case $Y$ is a $\PP ^{1}$ bundle over an elliptic curve 
or over $\PP^1$, one can easily construct a dominating map from 
$\PP^1\times\CC ^1$ and hence from $\CC^2$ which respects the bundle structure
(even algebraically in the latter case).   
In fact, by composing with the map 
\begin{equation}\label{eq:C^2}
(\pi^1, h\pi^2):\CC^2\rightarrow\CC^2
\end{equation}
where $h:\CC\rightarrow\CC$ is holomorphic with prescribed zeros
(which we can do by Weierstrass' theorem)
and $\pi^1,\ \pi^2$ are the respective projections, we can
arrange to have the dominating map miss any finite subset in
$Y$. Choosing this finite subset to be the set of indeterminacies of
the birational map from $Y$ to $X$, this dominating map lifts to give
a dominating map into $X$.  Since $\PP^1$ admits no holomorphic
differentials and is simply connected, we obtain, respectively,
$$
q(X)=q(Y)=q(C)=g \ \ \mbox{and}\ \ \pi_1(X)=\pi_1(Y)=\pi_1(C).
$$
Coupling this with the fact that the
fundamental group of a curve of genus greater than 1 is not a 
finite extension of an abelian group
gives us the following.

\begin{theorem} If $X$ is a projective surface with
${\kappa}(X) = - \infty$, then the following are equivalent.

\begin{itemize}
\item[\rm (a)] $X$ is dominable by $\CC ^ {2}$.
\item[\rm (b)] $q(X): = h ^ {0} (\Omega _{X} ) < 2$.
\item[\rm (c)] $X$ admits a Zariski dense holomorphic image of $\CC$.
\item[\rm (d)] $\pi_1(X)$ is a finite extension of an abelian group.
\end{itemize}
\end{theorem}


\subsection{Elliptic fibrations}  \label{subsec:elliptic}

If $X$ is any compact non-projective surface with $a(X)\neq 0$, then 
$X$ is an elliptic surface by \cite{Ko2}. 
Also, if $X$ is projective and $\kappa(X) = 1$, then $X$ is again an
elliptic surface by classification.  Hence the only remaining cases of
$\kappa \neq 0$ and $a \neq 0$ are elliptic surfaces.  In this section
we resolve completely the case of elliptic surfaces.

\begin{definition} An {\it elliptic fibration } is a proper
holomorphic map from a surface to a curve whose general fiber is an
elliptic curve, i.e., a curve of genus one. Such a surface is called
an elliptic surface. 
An elliptic fibration is called {\it relatively minimal} if
there are no $(-1)$-curves on any fiber.
\end{definition}

Note that an elliptic fibration structure on a minimal surface must be
relatively minimal.\medskip

Let  $f :  X \rightarrow C$ be a
fibration (i.e. a proper holomorphic map with connected fibers) 
between complex manifolds $X$ and
$C$. If $f':X'\ra C'$ is another map where $C'\subseteq C$, then
a map $h:X'\ra X$ is called fiber-preserving if $f\circ h=f'$. If
rank$(d f ) = \mbox{dim} \ C$ at every point on a fiber $X _{s} =
f ^ {-1} (s)$, then $X _{s}$ must be smooth by the implicit function
theorem.  If rank $(d f ) < \mbox{dim} \ C$ somewhere on $X _{s}$,
then $X _{s}$ is called a singular fiber.
Outside the singular fibers,
all fibers are diffeomorphic by Ehresmann's theorem.

\smallskip

In the case $f$ is a fibration of a surface $X$ over a curve $C$, then
each fiber, as a subscheme via the structure sheaf from $C$,
is naturally  an effective divisor on $X$ as follows.
We write $X _{s} = \sum n _{i} C _{i}$ 
where each $C _{i}$ is the $i$-th component of the fiber 
$(X_{s})_{red}$ (without the scheme structure) and where
$n _{i} -1$ is the vanishing order of $d f$ for a generic point
on $C_{i}$.
The positive integer coefficient $n_i$ is called the multiplicity
of the $i$-th component.
The multiplicity of a fiber $X _{s} = \sum n _{i} C _{i}$ is defined
as the greatest common divisor $n_s$ of $\{ n _{i} \}$.  A fiber 
$X_{s}$ with $n_s>1$ is called a multiple fiber.  A smooth fiber is 
then a fiber of 
multiplicity one having only one component.  The singular fibers form a 
discrete set in $X$ by analyticity. We will assume this setup for 
$X$ and $C$ from now on.

Let $\alpha:\tilde{C}\ra C$ be 
a finite proper morphism. The ramification index
at a point $\tilde s\in \tilde C$ is defined as the vanishing order of 
$d\alpha$ at $\tilde s$ plus one. Suppose $\alpha$ has ramification index
$n_s$ at every point above $s\in C$ and suppose that this 
is true for every $s\in C$.  Then,
according to  \cite[III, Theorem 9.1]{BPV}, pulling back the
fibration via this ramified cover yields an 
unramified covering $\tilde{X}$ over
$X$. Also, the resulting
fibration $\tilde X\ra \tilde{C}$ no longer has any multiple fibers.  
Such a ramified covering $\tilde{C}$ is called an orbifold covering
of $C$ with the given branched (orbifold) structure on $C$. 
More generally we have:

\begin{definition} \label{:D1} Given a curve $C$ with an assignment of a
positive integer $n_s$ for each $s\in C$ such that the set 
$S=\{s\in C|\ n_s>1\}$ is discrete in $C$, define 
$D = \displaystyle{\sum _{n _{s} >1}} 
\left( 1 - \frac{1}{n _{s}} \right) s$.
Suppose $\alpha:\tilde{C}\ra C$ is a holomorphic map such that
$\alpha:\tilde C\setminus \alpha^{-1}(S)\ra C\setminus S$ is an unramified
covering and such that, for each point $s\in S$,  
every point on $\tilde{C}$ above $s$ has ramification index
$n_{s}$. Then $\tilde C$ is called an orbifold covering of
the orbifold $(C,D)$.  If also $\tilde C$ is simply connected, 
then $\tilde{C}$ is called a uniformizing orbifold covering.
A fibration over $C$ defines a natural (branched) orbifold structure $D$
on $C$ by assigning $n_s$ to be the multiplicity of the fiber at $s$ of the
fibration.
\end{definition}

Therefore, we have the following:

\begin{prop}  \label{prop:L1} Let $X$ be a fibration over $C$.
Let $n_s$ denote the multiplicity of the fiber
$X_s$ for every point $s\in C$, thus endowing $C$ with an orbifold
structure $D$ as above. Let $\tilde{C}$ be an orbifold covering of $(C,D)$.
Then the pull back fibration $\tilde{X} \rightarrow \tilde{C}$ has no
multiple fibers and $\tilde{X} \rightarrow X$ is an
unramified holomorphic covering map.
\end{prop}

\medskip

\subsubsection{The Jacobian Fibration} 

We first begin with a preliminary discussion in the absolute case, the
case where the base is just one point.

\smallskip

Let $Z$ be a one dimensional subscheme (or a curve)
in a complex projective surface. The arithmetic genus of $Z$, defined by
$p_a(Z)= h ^ {1}  ( {\cal O} _{Z})
:= \mbox{dim} _{\CC} H ^ {1} ( {\cal O} _{Z} )$, 
is equal to the geometric genus
when $Z$ is smooth. Assume now that $Z$ is an arbitrary fiber in an
elliptic fibration. Since $p_a$ is an invariant in any algebraic family 
of curves (\cite[III, cor. 9.13]{Ha}), we have $p_a(Z)=1$ and so 
$H^1( {\cal O} _{Z} )=\CC.$ From the exponential exact sequence 
$0 \rightarrow \ZZ \rightarrow {\cal O}\rightarrow {\cal O} ^ {\ast} 
\rightarrow 0$, 
we construct the cohomology long exact sequence over $Z$ to deduce:
$$
\begin{array}{clclclclc}
0  &\rightarrow & H ^{1 } (Z, \ZZ ) & \stackrel{i}{\rightarrow} & 
H ^{1 } ( {\cal O} _{Z}
) & \rightarrow & H ^ {1 } ( {\cal O} _{Z} ^ {\ast} ) &
\stackrel{\delta}{\rightarrow} & H ^ {2} (Z, \ZZ ) \rightarrow 0 \\
   &\ & || & & || & & || & & || \\
   &\ &\mbox{a } \ZZ\mbox{-module} & & \CC & & \mbox{Pic} (Z) & & \ZZ
\end{array}
$$
{\bf Fact:} (Let $Z$ be non-singular.)
Pic$(Z)$ is naturally identified with the space of holomorphic line bundles 
over $Z$, which, in our case of $p_a=1$, is a 1-dimensional complex Lie group 
under tensor product. Every line bundle $L$ can 
be written as ${\cal O}(E)$ for some divisor $E=\sum a_i s_i$ 
($a_i\in \ZZ, s_i\in Z$) and $\delta(L)=\deg E:= \sum a_i$.

\begin{definition}
$\mbox{\rm Pic} ^ {0} (Z):= \ker \delta$ is the subgroup of {\rm Pic}$(Z)$ of
line bundles $L$ with trivial first Chern class $c _{1}(L):=\deg(L)$.
\end{definition}

If $Z$ is a smooth elliptic curve with a base point $\sigma$, we can construct 
a group homomorphism from $Z$ to Pic$^0(Z)$ by the map
$$
x \in Z \stackrel{f}{\mapsto} {\cal O} (x - \sigma ) 
\in \mbox{Pic} ^ {0} (Z).
$$

\begin{lemma}  \label{lemma:L2}
The map $f$ is holomorphic, one-to-one and hence onto. 
\end{lemma}

\noindent
{\bf Proof:} As $f$ is holomorphic by construction, we need to
prove only that it is one-to-one.
Assume not, so that ${\cal O} (x - \sigma ) = {\cal O} ( x ^ {\prime} -
\sigma)$ where $x \neq x ^{\prime}$. Then 
${\cal O} (x- x ^{\prime} )$ corresponds to the
trivial line bundle over $Z$ and so $Z$ has a rational 
function with a simple pole at $x ^ {\prime}$ and a simple zero at $x$. 
This gives a 1-1 and hence surjective holomorphic map
from $Z$, which has genus $1$, to $\PP ^ {1 }$, which has genus $0$. 
This is a contradiction.  $\BOX$

\noindent
{\bf Note:} $\mbox{Pic} ^ {0} (Z) = 
H ^{1 } ( {\cal O} _{Z} ) / i(H ^{1 }( Z, \ZZ )).$

\bigskip

We now return to the case in which the base is a curve.

\smallskip

Given an elliptic fibration $f:  X \rightarrow C$
without multiple fibers, one can construct
a relative version of $\mbox{Pic} ^ {0}$ as follows 
(see \cite[p. 153]{BPV}). We first form the  ${\cal O} _{C}$ module 
$${\cal J}\! ac ( f ) = f _{\ast 1} ( {\cal O} _{X} ) /
f _{\ast 1} \ZZ$$ over $C$.
Since $p_a( X _{s} ) = 1$
for every fiber, it follows that $f _{\ast 1} ( {\cal O}_{X} )$ 
is locally free of rank $1$ (by a well known theorem of Grauert)
and hence is the sheaf of sections of a line bundle $L$ over ${C}$.
Hence ${\cal J} \! ac ( f )$ corresponds
to the sheaf of sections of 
	$$\mbox{Jac}( f ):= L / f _{\ast 1} \ZZ,$$ 
which is a holomorphic fibration of complex Lie
groups with a zero section (see \cite{Ko2}, compare also \cite[V.9]{BPV}).  
Note that when $X _{s}$ is smooth elliptic, $( f _{\ast 1} \ZZ ) _{s} =
H ^{1 } ( X _{s} , \ZZ )$ which embeds in 
$L_s=H^1({\cal O}_{X_s})=\CC$. 
So Jac$( f ) _{s} = \mbox{Pic} ^ {0}
(X _{s} )$.  Note also that Jac$( f )$ is a holomorphic quotient of  
a line bundle $L$ over $C$.

\smallskip

We have the following theorem from Kodaira \cite{Ko2} (see
\cite[V9.1]{BPV}).

\begin{prop}   \label{prop:L3}
Let $f : X \rightarrow C$ be a relatively minimal elliptic fibration
over a curve $C$ with a holomorphic
section $\sigma : \ C \rightarrow X$.  Let $X
_{\sigma} ^ {\prime}$ consist of all irreducible components of fibers
$X_s$ not meeting $\sigma (C)$, and let $X^{\sigma} = X \setminus X
_{\sigma} ^ {\prime}$.  Then there is a canonical fiber-preserving
isomorphism $h$ from {\rm Jac}$( f )$ onto $X^{\sigma}$ mapping the
zero-section in {\rm Jac}$( f )$ onto $\sigma (C)$.
\end{prop}
\medskip 

Hence it is useful to construct holomorphic sections of elliptic fibrations
for which we develop the following key lemma.

\begin{lemma} \label{lem:sigma} Given a relatively minimal elliptic fibration 
$f: X\ra C$ without multiple fibers, assume $C$ is non-compact.
Then $f$ has a holomorphic section. Furthermore, given a countable
subset $T$ of $X$ whose image $f(T)$ is discrete in $C$, 
the section can be chosen to avoid $T$.
\end{lemma}

\noindent{\bf Proof:} From  Kodaira's table of non-multiple
singular fibers (\cite{Ko2} or \cite[Table 3 p. 150]{BPV}), we see that
every fiber which is not multiple in a relatively minimal elliptic fibration
has a component of multiplicity one.  So, every point on
${C}$ admits a neighborhood with a section.
We now choose a locally finite good covering of $C$ by open sets $U _{1} ,
U _{2} ,...,$ with sections $\tau _{1} , \tau _{2} ,...$ of 
${f} |_{U _{1}} , {f} |_ {U _{2}} ,...,$ respectively.  We may further
stipulate that there are no singular fibers on the intersection of any two
$U_j$'s.

Let $ {L} = {f} _{\ast 1} {\cal O} _{X}$, which is a holomorphically trivial
line bundle over $C$ since $C$ is Stein.  
Let $U \subseteq C$ be open and $\tau ^{\prime}
\in H ^{0} (U, {L} )$ a section. If $\tau$ is
a section of ${f} |_U$, then  we can
form the section $\tau+\tau ^{\prime}$ of ${f} |_U$
by proposition~\ref{prop:L3}.  By the same
proposition and the fact that all fibers are elliptic curves over $U
_{i} \cap U _{j}$, there is a section 
$\tau _{ij} ^ {\prime} \in H ^{0} (U _{i} \cap U _{j} , {L} )$ such that
$\tau _{i} + \tau _{ij} ^ {\prime} = \tau _{j}$ on $U _{i} \cap U _{j}$.

As $\{\tau_{ij}^{\prime}\}$ satisfies the cocycle condition, so does
$\{-\tau_{ij}^{\prime}\}$.
By the solution to the classical additive Cousin problem (or from the fact that
$H ^{1 } ( \{ U _{i} \} ,L) = H ^{1 } ( C ,L) = H^1(C,{\cal O}) =0$ by
Leray's theorem, Dolbeault's isomorphism, and the fact that
$C$ is Stein), one can find
holomorphic sections $\tau _{i} ^ 
{\prime} \in H ^{0} (U _{i} , {L} )$ such that $\tau _{i} ^ {\prime } -
\tau _{ij} ^ {\prime} = \tau _{j} ^ {\prime}$.  Then $\tau _{i} + \tau _{i} ^
{\prime} = \tau _{j} + \tau _{j} ^ {\prime}$ on $U_i\cap U_j$ for all $i,j$.
This gives rise to a global section of $f : X \rightarrow C$.

Given such a global section,
proposition~\ref{prop:L3} gives a fiber-preserving
dominating map $F: L\ra X$ where $F^{-1}(x) \subset L$ is at most a countable
discrete set for all $x$ in $X$. Hence $F^{-1}(T)$
is also a countable set and is supported on the fibers of $L$ over $f(T)$.
For each $s\in f(T)$, therefore, we may choose a point $q_s$ in
$L\setminus T$.  As $L$ is isomorphic to the trivial 
line bundle ($C$ being non-compact), 
the classical interpolation theorems of Mittag-Leffler
and of Weierstrass give us a holomorphic section $\sigma$ of $L$
with the prescribed value $q_s$ for all $s\in T$. But then
$F\circ \sigma$ is a section of $f$ which avoids $T$.
This completes the proof.  $\BOX$

\subsubsection{Theorem \ref{thm:1} in the case of elliptic fibrations}

\begin{theorem} \label{thm:elliptic-fibration}
Let $f : X \rightarrow C$ be a relatively minimal 
elliptic fibration with a finite number of multiple fibers. Assume that
$C$ is a Zariski open subset of 
a projective curve $\bar C$.
Let $n _{s}$ be the
multiplicity of the fiber $X _{s}$.
Then the following are equivalent.

\begin{itemize}
\item[\rm (a)] $X$ is dominable by $\CC ^ {2}$.
\item[\rm (b)] $\chi: = 2 - 2 g( \overline{C} ) - \#
( \overline{C} \setminus C) - \displaystyle{\sum _{n _{s} \geq 2}}
(1 - \frac{1}{n _{s}}) \geq 0$.
\item[\rm (c)]  There exists a holomorphic map of $\CC$ to $X$ whose
image is Zariski dense. 
\end{itemize}
\end{theorem}

\noindent
{\bf Remark 1:}  $\chi=\chi(C,D)$ 
is the orbifold Euler characteristic of $(C,D)$.
It can be written as 
$\chi (C,D)= 2 - 2g ( \overline{C} ) - \displaystyle{\sum _{s \in
\overline{C} }} (1 - \frac{1}{n _{s}})$ 
if we set $n _{s} = \infty$ for $s \in \overline{C} \setminus C$
(where  $\frac{1}{\infty}=0$).  Hence, if we complete the $\QQ$-divisor
$ D= \displaystyle{\sum _{s \in C}} ( 1 - \frac{1}{n _{s}} ) s\,$
to 
$\overline{D} = \displaystyle{\sum _{s \in \overline{C} }} (1 -
\frac{1}{n _{s}} ) s$ on $\overline{C}$, then
$$ \chi(C,D)= 2 - 2g ( \overline{C} ) - \deg \overline{D}.$$

\medskip \noindent
{\bf Proof of theorem:}  The pair $(C,D)$ 
defines an orbifold as given in definition~\ref{:D1}. We will show
that (a) holds if $\chi (C,D)\geq 0$ while property C
holds for $X$ (that is, (c) fails to hold) if $\chi (C,D)<0$. This
will conclude the proof.

{}	From the classical uniformization theorem
for orbifold Riemann surfaces (see, for example, 
\cite[IV 9.12]{FK}), $(C,D)$ has
a uniformizing orbifold covering $\tilde C$ which is $\PP^1$, $\CC$ or $\DD$ 
according to $\chi (C,D) > 0$, $\chi (C,D)=0$ or $\chi (C,D) < 0$
respectively, unless $\overline{C} = \PP ^ {1 }$ and $\overline{D}$
has one or two components.  In the latter (``unless'') case, we simply
redefine $C$ to be the complement of the components of $\overline{D}$
in  $\PP ^{1 }$ and reset $D$ to be $0$, shrinking $X$ as
a result. We can do this because it does not change the fact that
$\chi(C,D) \geq 0$ and because
once we show that the resulting $X$ is dominable
by $\CC^2$, the original $X$ is also.

By pulling back the fibration to $\tilde{C}$, 
we obtain a relatively minimal elliptic fibration $Y$ over $\tilde{C}$.
Now, proposition~\ref{prop:L1} 
implies that the natural map from $Y$ to $X$ is an
unramified covering. Hence any holomorphic map from $\CC$ to $X$ 
must lift to a holomorphic map to $Y$. It follows that if $\tilde{C}=\DD$,
then any such map must lift to a fiber and hence its image in $X$ must
lie in a fiber. So, property C holds
and $X$ cannot be dominated by $\CC^2$ in this case.
 
	It remains to show that $X$ is
dominated by $\CC ^ {2}$ in the case $\tilde{C} = \CC$ or $\PP^1$
to complete the proof of this theorem.  Note that 
the latter case can be reduced
to the former by simply removing a point from $\tilde{C}$. Hence, we may take 
$\tilde{C}$ to be $\CC$ which is non-compact. Lemma~\ref{lem:sigma} 
now applies to give a section of the pullback fibration
$\tilde f : Y\rightarrow \CC$.
By proposition~\ref{prop:L3}, $Y$ is dominated by Jac$(\tilde f)$
which in turn is dominated by a line bundle $L$ over $\CC$ by construction.
Hence $X$ is dominated by $L=\CC^2$ 
(since any line bundle over $\CC$ is holomorphically
trivial) as required.  $\BOX$

Now, let $f':X'\ra C$ be an arbitrary elliptic
fibration. By contracting the $(-1)$-curves on the fiber,
we get a bimeromorphic map $\alpha$ from $X'$ to a
surface $X$ having a relatively minimal elliptic fibration 
structure over $C$. As before, $X$ defines an orbifold structure $D$ 
on $C$. If $X$ has an infinite number of multiple fibers or
if $C$ is not quasi-projective, then $\DD$ is the universal covering
of $(C,D)$ and conditions (a) and (c) of this theorem both fail for $X$.
Otherwise the above theorem can be applied to conclude that conditions
(a) and (c) are still equivalent for $X$.
Let $T$ be the indeterminacy set of $\alpha^{-1}$.
By examining the last paragraph of the
above proof, we see that lemma~\ref{lem:sigma}
actually applies to give us a dominating map from the
trivial line bundle $L$ over
$\tilde C$ to $X$, and the zero-section of $L$ maps to a section of $f$ that
avoids $T$. Composing with a self-map of $L$ given by a section of $L$
with prescribed zeros (just as in equation~\ref{eq:C^2}) then gives us a 
dominating map from $L$ to $X$ which avoids $T$. Hence, if $X$ is 
dominable by $\CC^2$, then $X'$ is also.
It is clear that $X'$ satisfies property C if $X$ does.
Hence, we obtain the following, which covers theorem~\ref{thm:1} in
the case of elliptic surfaces.

\begin{theorem} \label{thm:elliptic}
Let $f : X \rightarrow C$ be an 
elliptic fibration. 
Then conditions (a) and (c) of theorem 
\ref{thm:elliptic-fibration} above are equivalent for $X$; that is,
dominability by $\CC^2$ is equivalent to having a Zariski-dense holomorphic
image of $\CC$.
\end{theorem}

\vspace{-0.1in} Note that we do not require $C$ to be quasi-projective in
this theorem.

\medskip

\subsubsection{An algebro-geometric characterization}

In this section, we will give, without proof, a characterization of
dominability by $\CC^2$ for a projective elliptic fibration in terms of
familiar quantities in algebraic geometry and not involving the 
fundamental group.  Unfortunately, the condition given is not
straightforward nor does it seem very tractable. Hence, we will
leave the proof (which is based on the simple fact that the saturation
of the cotangent sheaf of the base,
pulled back by the fibration map, includes the orbifold cotangent sheaf
as a $\QQ$ subsheaf) to the reader. We will deal only with the case
of $\kappa=1$ since the other possibility of $\kappa=0$ contains the,
so far, problematic K3 surfaces. However, all surfaces with $\kappa =
0$ other than the K3's are dominable by $\CC^2$.
We note that from the classification list in section 2, a surface with
$\kappa=1$ is necessarily elliptic.
Before the statement of the following proposition, recall that a
vector sheaf is called big if it contains an ample subsheaf. Recall
also that a divisor in a surface is nef if its
intersection with any effective divisor is non-negative.

\begin{prop}
Let $X$ be a projective surface with $\kappa(X)=1$. Then $X$ is dominable by
$\CC^2$ if and only if there exists a nef and big divisor $H$ such 
that, for every nef divisor $N$ with $K_X N=0$, there exists a positive
integer $m$ with ${\rm S}^m\Omega_X(H-N)$ big.
\end{prop}

It is not difficult to extract a birational invariant out of this from
$\QQ$-subsheaves of the cotangent sheaf of such an elliptic surface; 
again we leave this to the interested reader.

In the remainder of this section, we give a 
more satisfactory and elementary characterization of
dominability, now in terms of the fundamental group.

\subsection{The fundamental group of an elliptic fibration}

	We begin with the remark that, except for our narrow focus on
holomorphic geometry, most of the results we obtain in this section 
are not presumed to be new.

	Let $ f: X\ra C$  be an elliptic fibration.
Then the fibration determines a branched orbifold structure $D$ on $ C$
as given in definition~\ref{:D1}.
Let $C^\circ$ be the complement of the set of branch points in $ C$. Then
$X^\circ= f^{-1}(C)$ is an elliptic fibration 
defined by $f^\circ= f|_{X^\circ}$, which has no multiple fiber. Let
$X'$ be the complement of the singular fibers in $\bar X$.  
Then $f'=f|_{X'}$ defines a smooth fibration over a curve $C'\subseteq C$,
and is therefore differentiably locally trivial by Ehresmann's theorem.
We have the following commutative diagram.
\begin{equation} \label{diag:ell}
\begin{array}{rcrcl}
X' & \hookrightarrow & X^\circ  & \hookrightarrow &  X \ \\
f' \downarrow \, & &  f^\circ \downarrow \, & & \,\downarrow   f \\
C' & \hookrightarrow & C^\circ  & \hookrightarrow &  C  
\end{array}
\end{equation}

We first observe the following trivial lemma
for our consideration of $\pi_1(X)$. Throughout
this section, all paths are assumed
to be continuous.

\begin{lemma}\label{lem:cod2} Assume that we are
given a real codimension two subset $W$ of $X$ and a path 
$\nu:[0,1]\ra X$ such that $\nu(0)$ and $\nu(1)$ lies outside $W$.
Then $\nu$ is homotopic to a path that avoids $W$
keeping the end points fixed.
\end{lemma}

\noindent{\bf Proof:} We first impose a metric on $X$. Since $[0,1]$ 
is compact,
there is an integer $n$ such that $\nu([(i-1)/n,i/n])$ is 
contained in a geodesically convex open ball $B_i$ for all 
$i\in \{1,2,...,n\}$.
Then the intersection of these balls are also  geodesically convex and,
in particular, connected. Now replace $\nu(i/n)$ by a point in 
$B_i\bigcap B_{i+1}\setminus W$ for each integer $i\in [1,n-1]$. Then replace 
$\nu|_{[(i-1)/n,i/n]}$ by a path in $B_i\setminus W$ connecting $\nu((i-1)/n)$
with $\nu(i/n)$, for each integer $i\in [1,n]$. This is possible because
the complement of $W$ in each of the open balls is connected as $W$ is of real
codimension two in them. Since the balls are contractible and intersect
in connected open sets, we see that the new path is homotopic to the original
one fixing the end points but now avoids $W$. \BOX

If the path $\nu$ given above has the same end points, that is $\nu(0)=\nu(1)$,
then we call $\nu$ a loop. We will often identify $\nu$ with its image.

\medskip

For the next two propositions, we observe from Kodaira's table of 
singular fibers (see \cite[V.7]{BPV}) that, for a fiber $X_s$ of an
elliptic fibration (as  a topological space or a simplicial complex),
$\pi_1(X_s)$ is either $\ZZ\oplus\ZZ$ (corresponding
to a nonsingular elliptic curve), $\ZZ$ (corresponding to the (semi-)stable
singular fibers), or the trivial group (corresponding to the 
other singular fibers).

\begin{prop} \label{prop:elliptic-pi}
Let $f: X\ra C$ be an elliptic fibration. In the case $C=\PP^1$, let
$X_\infty$ be a multiple fiber if one exists. Assume
$f$ has no multiple fibers except possibly for $X_\infty$
and that $C$ is simply connected.
Then $\pi_1(X)$ is a quotient of $\pi_1(X_s)$ for every fiber outside 
$X_\infty$. In particular,  $\pi_1(X)$ is abelian.
\end{prop}

\noindent
{\bf Proof:} Since contracting $(-1)$-curves does not change the fundamental
group, we may assume without loss of generality that $f$ is relatively 
minimal. Let $X_s$ be an arbitrary fiber. 
Being a CW-subcomplex of $X$, it
is a deformation retract of a small neighborhood $U$  which we may assume 
to contain a smooth fiber $X_{s'}$ nearby.
Since $X$ is path connected, we can choose any base point
in considering its fundamental group. Fix then a base point 
$q\in X_{s'}$ and a loop $Q$ with this base point.
We will show that $Q$ is pointed homotopy equivalent in $X$ 
to a loop in $X_{s'}\subset U$. The theorem then follows as $X_s$ is
a deformation retract of $U$.

Since the singular fibers form a real
codimension two subset, we can modify $Q$ to avoid
them up to pointed homotopy equivalence by lemma \ref{lem:cod2} above.
In the case $C=\PP^1$ but $X_\infty$ is not already given, 
let $X_\infty$ be a fiber outside $U$ and 
$Q$. Since  every homotopy (of $Q$) in $X\setminus X_\infty$ is also one
in $X$, we may safely replace $X$ by $X\setminus X_\infty$ so that $C$
becomes contractible in this case. Hence, we may assume in all cases
that $C$ is contractible and that $Q$ is a loop in $X'$, the
complement of the singular fibers in $X$. So lemma 
\ref{lem:sigma} and proposition \ref{prop:L3} apply to give an isomorphism
from Jac$(f)$ to $X$ with parts of the singular fibers complemented. Hence,
we get, by construction of Jac$(f)$, a map $\theta$ from a holomorphically
trivial line bundle $L$ over $C$ to $X$ which is an unramified covering above
$X'\subset X$. Fixing a point $q_0\in \theta^{-1}(q)\subset L_{s'}$, we see that
$Q$ can be lifted to a path $\tilde Q$
in $L$ from $q_0$ to a point $q_1\in L_{s'}$
by the theory of covering spaces. As $C$ is contractible, there is a 
homotopy retraction of $L$ to $L_{s'}$ which provides a pointed homotopy of
$\tilde Q$ to a path in $L_{s'}$. Pushing down this homotopy (via $\theta$)
to $X$ gives a pointed homotopy from $Q$ to a loop in $X_{s'}$ as required.
\BOX

Looking back at the above proof, we see that we can reach the same 
conclusion by allowing $X_s$ to be a multiple fiber as long as $C$
is contractible and $X$ is free of other multiple fibers. This can
be done by contracting the loop $Q$, as given in the proof, but only to the
neighborhood $U$ of $X_s$ before homotoping to $X_s$ via the deformation
retraction of $U$ to $X_s$. 
Of course, $X_s$ as stated in the theorem is no longer
arbitrary in this case as it is a multiple fiber. If $C=\PP^1$ and
$D$ has two components (corresponding to $X$ having two multiple fibers),
we can remove one of the components (corresponding to removing one multiple
fiber from $X$) for the same conclusion. We recall that 
$$D=\sum_{s\in C}\Big(1-\frac{1}{n_s}\Big)\,s$$ 
defines the orbifold structure on 
$X$ where $n_s$ is the multiplicity of the fiber $X_s$.
Hence, we get a complement
to the above proposition.

\begin{prop}\label{prop:star} Let $f:X\ra C$ be an elliptic fibration 
defining the orbifold structure $D$ on $C$.
If $C=\PP^1$ and $D$ has one or two components,
or if $C$ is contractible and $D$ has one component, then $\pi_1(X)$ is
a quotient of $\pi_1(X_s)$ for every component $s$ of $D$. Hence, $\pi_1(X)$
is abelian in these cases.\end{prop}\BOX

\subsubsection{Monodromy action as conjugation in the fundamental group}
Although it is not absolutely necessary, some familiarity 
with the notion of monodromy and vanishing
cycles used in geometry may be useful for reading this section.

Let the setup be as in diagram \ref{diag:ell} and let $X_r$ be a 
non-singular fiber. Fix 
a base point $q$ in $X_r$ for all fundamental group considerations
from now on. 
There is an action of $\pi_1(X' ,q)$ on $\pi_1(X_r,q)$ via
the monodromy action
which, in the case $C'$ is not $\PP^1$, is just the 
conjugation action in 
$\pi_1(X' ,q)$. Indeed, in this case, we have the following exact
sequence from the theory of fiber bundles (or from elementary covering 
space theory) 
\begin{equation}\label{seq:pi}
0\ra \pi_1(X_r,q) \ra \pi_1(X' ,q) \ra \pi_1(C',r)\ra 0,
\end{equation}
from which we deduce that the monodromy action is really an action of
$\pi_1(C',r)$ on $\pi_1(X_r,q)$ since the latter is abelian.  

In general, we will let $H$
denote the image of $\pi_1(X_r,q)$ in $ \pi_1(X ,q)$ under the inclusion
of $X_r$ in $X$. It is easy to see that $H$ is a normal subgroup in $\pi_1(X)$
(by the definition of the monodromy action). 
In this paper, we will be mainly interested in the monodromy action on $H$.
As opposed to the 
usual case of the monodromy action on the homology level, this
action need not be trivial, unless we know, for example, that $\pi_1(X, q)$
is abelian. Hence, it is of interest for us to know how far $\pi_1(X, q)$ is
from abelian.

With the same setup, suppose $(C,D)$ has a 
uniformizing orbifold cover $\tilde C$.
This is the case 
unless $\bar C=\PP^1$ and $\bar D$ has one or two components,
again by the uniformization theorem (\cite[IV 9.12]{FK}). In the latter
cases, proposition~\ref{prop:star} and proposition~\ref{prop:elliptic-pi}
tell us that $\pi_1(X)$ is abelian so that 
the monodromy action on $H$ is 
trivial. In all other cases, let $\tilde f:\tilde X\ra \tilde C$ 
be the pullback fibration. 
Proposition~\ref{prop:L1} implies that $\tilde X$ is 
an unramified cover over $X$. Let $R$ be the covering group and 
$G=\pi_1(X)$. From the theory of covering
spaces, we know that $G$ is an extension of 
$\pi_1(\tilde X)$ by $R$. Since  $\pi_1(X_r)$ surjects to
$\pi_1(\tilde X)\subset \pi_1(X)$ by proposition~\ref{prop:L1}, we see that
$$H=\pi_1(\tilde X).$$ Note that $R$ is
a quotient of $\pi_1(C^\circ)$
and hence also of $\pi_1(C')$, allowing us to identify the conjugation action 
of $R$ on $H$ with the monodromy action. Hence,
we have the following exact sequence 
(which we can regard as a quotient of the exact sequence \ref{seq:pi})
\begin{equation}\label{seq}
0\ra H\ra G \ra R \ra 0.
\end{equation}
The following proposition tells us that this monodromy action
on $H$ via loops in $X'$, which induces the conjugation action
of $R$ on $H$, depends only on 
the pointed homotopy class of the image of these
loops in $C$. Hence, the monodromy action on $H$ is really an action
by the group $\pi_1(C)$, which is a quotient of $R$.
In particular, it tells us that
the action is trivial when $C$ is simply connected. This is the closest
analogue, on the level of $\pi_1$, 
of the fact that vanishing cycles are vanishing on
the level of homology.

\begin{prop}\label{prop:homotopy} Let $f:X\ra C$ be an elliptic fibration.
Let $X_r$ be a non-singular fiber with a base point $q$. If $\alpha$,
$\beta$ and $\gamma$ are loops based at $q$ with $\alpha$ in $X_r$, and
$f\circ \beta$ is pointed homotopic to $f\circ \gamma$ in $C$, then
$\beta^{-1}\alpha\beta$ is pointed homotopic in $X$ to 
$\gamma^{-1}\alpha\gamma$.\end{prop}

\noindent
{\bf Proof:} We may assume, via lemma \ref{lem:cod2}, that $\beta$ and
$\gamma$ lie in $X'$. Let $h:[0,1]\times [0,1]\ra C$ be a pointed 
homotopy between $f\circ \beta$ and $f\circ \gamma$, which exists by
assumption. Note that
$$(f\circ\beta)(f\circ\gamma)^{-1}=h\big(\partial ([0,1]\times [0,1])\big)$$ 
as loops up to pointed homotopy equivalence,
where $\partial$ means the oriented boundary. Our conclusion would follow
if we show that the monodromy action of this latter loop, call it
$\mu$, on $\alpha$ is trivial in $\pi_1(X)$.

By compactness of $h([0,1]\times [0,1])$, there is a partition
$\{0 = a_0 < a_1 < ... < a_n = 1\}$  of $[0,1]$ such that 
$h([a_{i-1},a_i]\times [a_{j-1}, a_j])$ is contained in an open disk $D_{ij}$
containing at most one branch point
and such that the loop
$$\mu_{ij}:=h\big(\partial ( [a_{i-1},a_i]\times [a_{j-1}, a_j])\big)$$
lies in $X'$, for all $i, j\in \{1,2,...,n\}$.
Since $\pi_1\big(f^{-1}(D_{ij})\big)$ 
is abelian by proposition \ref{prop:star},
the monodromy action of $\mu_{ij}$ on any pointed loop in the fiber
is trivial in $\pi_1(f^{-1}(D_{ij}))$, and hence in $\pi_1(X)$ as well,
for all  $i, j\in \{1,2,...,n\}$.
Our result now follows from the fact that the monodromy action of $\mu$
is just the sum of the monodromy action of the $\mu_{ij}$'s. \BOX

We can do a bit better when $X$ is compact.

\begin{lemma} \label{lemma:IST}
With the setup as in the above proposition, assume further that 
either $X$ is compact or $X$ is a holomorphic fiber bundle over $C$. 
Then an integer $m$ exists, independent of $\beta$, such
that $\beta^m$ commutes with $\alpha$ in $\pi_1(X)$.
\end{lemma}

\noindent
{\bf Proof:} Let $H$ be the image of $\pi_1(X_r)$ in $\pi_1(X)$. 
We may assume, as before, that $f$ is relatively minimal. 

If $X$ has a singular fiber, then $H$ is cyclic and hence the result 
follows from the fact that the automorphism group of a cyclic group
is finite.

If $X$ is a holomorphic fiber bundle over $C$, then the monodromy actions can
be realized as holomorphic automorphisms of the fiber. The group of 
such automorphisms is a finite cyclic extension of the group of lattice
translations (this can be deduced easily or determined from the table
in V.5 of \cite{BPV} listing such groups). Hence every monodromy
action up to a power is a translation on the fiber, which therefore
leaves every element of $\pi_1(X_r)$ invariant. 

If $X$ is compact
and has no singular fibers, then it is a holomorphic fiber bundle by
Kodaira's theory of Jacobian fibrations. So the result follows by
the last paragraph. \BOX

If $X$ is non-compact and $f$ is
algebraic without singular fibers, then the 
conclusion of this lemma may no longer hold.
Neverthless, we can embed $X$ in a projective 
surface $\bar X$, which is again elliptic.
Deligne's Invariant Subspace Theorem  \cite{De} 
implies that elements in $\pi_1(X_r)$ which vanish
in $\pi_1(X)$ are generated 
over $\QQ$ by commutators of the form given by this lemma. But we can
deduce this directly from the fact that the 
abelianization of $\pi_1(\bar X)$ must have even rank so that either $H$
lies in the center of $\pi_1(\bar X)$ (in the case when $\bar X$ is birational 
to an elliptic fiber bundle) or the commutator subgroup of 
$\pi_1(\bar X)$ generates $H$ over $\QQ$. In fact, Kodaira's theory
allows us to deduce a strong version of the Invariant Subspace Theorem
(in the case of elliptic fibrations)
which is valid even outside the algebraic category:

\begin{prop}\label{prop:kod} Let $f:X\ra C$ be an elliptic fibration
without singular fibers and such that $C$ is the complement of
a discrete set in a quasi-projective curve. Let $X_r$
be a fiber. Then either $f$ is holomorphically locally trivial or
$\pi_1(X_r)\otimes \QQ =H_1(X_r)\otimes \QQ$ is generated by the vanishing
cycles --- that is, by loops of the form $\alpha^{-1}\beta^{-1}\alpha\beta$
(naturally identified as elements of $\pi_1(X_r)$ via monodromy)
in the notation of proposition~\ref{prop:homotopy}, where $\alpha$
is a loop in $X_r$.
\end{prop}

We remark that a weaker form of this proposition is 
in fact due to Kodaira and is disguised in the proof of
theorem 11.7 in \cite{Ko2}. We will follow his method, almost verbatim,
in our proof. 

\medskip

\noindent
{\bf Proof:}  We begin with some preliminaries concerning the 
period function $z(s)$, which takes values in the upper half plane.
Recall that, as far as monodromy actions are concerned, we can identify
$\beta\in \pi_1(X)$ with an element of $\pi_1(C)$, which we will denote
again as $\beta$ by abuse of notation.
 
By theorem 7.1 and theorem 7.2 of \cite{Ko2} (neither of
which requires the additional
assumption of that section concerning the compactification), we have
a multivalued holomorphic 
period function $z(s)$ on $C$ with positive imaginary
part such that, under the monodromy representation $(\beta)\in$SL$(2,\ZZ)$ of 
$\beta \in \pi_1(C)$ as an automorphism of
the lattice $H_1(X_r)$ with a fixed choice
of basis, $z(r)$ transforms as 
$$
\beta_*:z(r)\longmapsto \frac{az(r)+b}{cz(r)+d},\ \ \ \ \mbox{where}\ \ \ \
(\beta)=\left(
\begin{array}{cc}
a&b\\
c&d
\end{array}
\right)\in \mbox{SL}(2,\ZZ)
$$
under our choice of basis.
By definition, $(1, z(s))$ is the period defining the elliptic curve $X_s$
via analytic continuation of $(1, z(r))$, which is fixed by our choice of
basis on $H_1(X_r)$ (see equation 7.3 in \cite{Ko2}).

With a choice of basis over the point $r$ fixed,
we can regard the period function $z$ as a single valued holomorphic
function on the universal cover $\tilde C$. Also, we can naturally identify
$\pi_1(C)$ with the covering transformation group of $\tilde C$ over $C$.
Then we have (see equation 8.2 in \cite{Ko2})
$$
z(\beta(\xi))=\beta_*z(\xi)=\frac{az(\xi)+b}{cz(\xi)+d},
\ \ \ \ \mbox{where}\ \ \ \
(\beta)=\left(
\begin{array}{cc}
a&b\\
c&d
\end{array}
\right)\ \ \mbox{and}\ \ \xi\in \tilde C.
$$
\smallskip

Let $M$ denote the submodule of $\pi_1(X_r)=H_1(X_r)=\ZZ\oplus\ZZ$
generated by the vanishing cycles. After a suitable change of basis,
we may assume that $M=n\ZZ\oplus m\ZZ\subseteq \ZZ\oplus\ZZ$, where
$m$ and $n$ are integers. If $M$ does not generate $H_1(X_r)$ over $\QQ$,
then either $m$ or $n$ must vanish. If $m$ vanishes, then we must have
$$
(\beta)=\left(
\begin{array}{cc}
1&b_\beta\\
0&1
\end{array}
\right)\ \ \ (\mbox{for some}\ b_\beta\in \ZZ),
$$
and therefore $z(\beta(\xi))=z(\xi)+b_\beta$ for all $\beta \in \pi_1(C)$.
Since the imaginary part of $z(s)$ is positive, exp$[2\pi iz(s)]$
defines a single valued holomorphic function on $C$ with modulus less
than $1$.  Hence, it must extend to a bounded holomorphic function on
the compactification $\bar C$ of $C$ and therefore must be constant.
It follows that $z(s)$ is constant and so the fibration is locally
holomorphically trivial. If $n$ vanishes, then
$$
(\beta)=\left(
\begin{array}{cc}
1&0\\
c_\beta&1
\end{array}
\right)\ \ \ (\mbox{for some}\ c_\beta\in \ZZ),
$$
and therefore 
$$1/z(\beta(\xi))=1/z(\xi)+c_\beta.$$
Hence, considering exp$[-2\pi i/z(s)]$ instead of exp$[2\pi iz(s)]$ 
gives us the same conclusion. This ends our proof.
\BOX

\medskip

	In order to study $G=\pi_1(X)$, 
we need some information about its quotient $R$. This
is fortunately a classical subject that we now turn to.

\subsubsection{Fuchsian groups versus elementary groups}

 In this section, we will collect some
basic definitions and facts that we will need 
about Kleinian groups. We refer the reader to 
\cite[IV.5-IV.9]{FK} and \cite[I-V]{Mas} for more details.

\smallskip

Let $(C,D)$ be an orbifold
and let $\tilde C$
be its uniformizing orbifold covering with covering group $R$ which acts
holomorphically on $\tilde C$.
Since $\tilde C=\PP^1, \CC$ or $\DD$, all of which have natural embeddings
into $\PP^1$, $R$ can be identified as a subgroup of the
group $\MM$ of holomorphic automorphisms of $\PP^1$, the group of Mobius 
transformations. So identified, $R$ becomes a Kleinian group; that is, 
a subgroup of $\MM$ with a properly discontinuous action at some
point, and hence in some maximum open subset $\Omega$, of $\PP^1$.
The set of points $\Lambda=\PP^1\setminus \Omega$ 
where $R$ does not act properly discontinuously is called the limit set of $R$.

An elementary group is a Kleinian group $R$ with no more than two points
in its limit set. Such a group acts properly discontinuously on 
$\Delta \subset \PP^1$, where $\Delta$ is $\PP^1$, $\CC$ or $\CC^\ast$.

By a Fuchsian group, we mean a Kleinian group $R$ with with
a properly discontinuous action on some disk $\DD\subset \PP^1$
such that $\DD/R$ is quasi-projective; that is, $\DD$
is the uniformizing orbifold covering of an orbifold $(C, D)$
where $C=\DD/R$ is quasi-projective. If $D$ has finitely many components,
then $(X, D)$ is known as a finite marked Riemann surface and $R$ is
called basic. The limit set of a Fuchsian 
group necessarily contains the boundary
of $\DD$ (which characterizes Fuchsian groups of the first kind
in the literature). It follows that a Fuchsian group cannot be an
elementary group.  We can also see this directly as follows.

\begin{lemma}\label{lem:elem} 
An elementary Kleinian group is not a Fuchsian group.
\end{lemma}

{\bf Proof:} Let $R$ be an elementary
Kleinian group, then $R$ acts properly discontinuously
on $\Delta=\PP^1,\ \CC$ or $\CC^\ast$ as a subset of $\PP^1$. 
If $R$ also acts on a disk $\DD\subset \PP^1$, then the boundary of
this disk with at most two points removed is contained in $\Delta$.
Since $R$ is properly discontinuous on $\Delta$, and hence on this 
punctured boundary, $\DD/R$ is not quasi-projective. Hence $R$ is not
Fuchsian.\BOX

\smallskip

The following is a direct consequence of the uniformization theorem.

\begin{prop}\label{prop:unif}
Let $(C,D)$ be a uniformizable orbifold where $D$ has a finite
number of components. Let $R$ be the uniformizing orbifold covering group of
$(C, D)$ properly regarded as a Kleinian group. Then $\chi (C,D)<0$ if
and only if $R$ is a Fuchsian group while $\chi (C,D)\geq 0$ if and only
if $R$ is an elementary group.
\end{prop}

The reader is cautioned that lemma~\ref{lem:elem} is not a corollary of
this since the definition of an elementary group is more general than
that given in this proposition.

\medskip

Concerning $R$ as an abstract group,  proposition~\ref{prop:unif}
and the basic theory of elementary Kleinian
groups (see \cite[V.C and V.D]{Mas} or \cite[IV 9.5]{FK}) gives:

\begin{prop}\label{prop:A} With the same setup as proposition~\ref{prop:unif},
assume that $\chi (C, D)\geq 0$. Then there is a finite orbifold covering 
$\tilde C$ of $(C, D)$ such that $\tilde C=\PP^1$, $\CC^\ast$ or an elliptic
curve. In particular, $R$ is a finite extension of a free abelian group of
rank at most two.
\end{prop}

Quoting \cite[V.G.6]{Mas}, using lemma~\ref{lem:elem} and 
proposition~\ref{prop:A}, we have:

\begin{prop}\label{prop:AA}
Let $R$ be a Fuchsian group
as defined above. Then $R$ is not a finite extension of an
abelian group. Hence, $R$ is not isomorphic to an elementary group as
an abstract group.
\end{prop}

\subsubsection{The fundamental group characterization in
theorem~\ref{thm:2}}\label{sec:pi_1}

Before stating the main theorem of this section, we need the following
proposition from \cite[V.5]{BPV}. We first note from 
the same source that an elliptic fiber
bundle over an elliptic curve is called a primary Kodaira surface if it
is not K\"ahler. A non-trivial free quotient of such a surface 
by a finite group is called a
secondary Kodaira surface. The fundamental group of such a surface is
unfortunately not a finite extension of an abelian group, even though
the surface is $\CC^2$-dominable.

\begin{prop}\label{prop:Ebundle} An
elliptic fiber bundle over an elliptic curve is either a primary Kodaira
surface, or a free and finite 
quotient of a compact complex 2-dimensional torus.
\end{prop}

Armed with this, we are ready to tackle our second main theorem,
theorem~\ref{thm:2}, in the case of elliptic fibrations. We will
state a more general theorem:

\begin{theorem} \label{thm:E2} Let $f: X\ra C$ be an
elliptic fibration with $C$ quasi-projective. Assume that $X$ is not
bimeromorphic to a free and finite quotient of a primary Kodaira
surface. Then $X$ is dominable by $\CC^2$ if and only if $\pi_1(X)$
is a finite extension of an abelian group (of rank at most 4).
\end{theorem}

\noindent
{\bf Proof:} With the assumptions as in the theorem, we let
$G=\pi_1(X)$ as before. 
By the same argument as that for theorem~\ref{thm:elliptic},
we may assume, without loss of 
generality, that $X$ is relatively minimal by contracting
the $(-1)$-curves (as $G$ is unchanged in this process). 
If $X$ has an infinite number of multiple fibers, then
the orbifold $(C,D)$ is uniformized by $\DD$ and so $X$ is not dominable
by $\CC^2$.
Proposition~\ref{prop:AA}
tells us that $R$ is not isomorphic to a finite extension of an abelian
group in this case. Hence, we may also assume that $D$ has only a finite number
of components for the rest of the proof. Theorem~\ref{thm:elliptic-fibration}
then applies and so it is sufficient to show that $\chi (C,D)\geq 0$
if and only if  $G$ is a finite extension
of an abelian group.

Assume that $\chi (C,D)<0$. If $(C, D)$ projectivizes to 
$(\PP^1, \bar D)$ (see the first
remark after theorem~\ref{thm:elliptic-fibration}
for the definition of $\bar D$),
then $\bar D$ must have more than two components by the definition
of $\chi$. Hence $(C,D)$ is uniformizable and 
we may apply proposition~\ref{prop:unif} and proposition~\ref{prop:AA} to
conclude that the orbifold uniformizing group $R$ of $(C,D)$ is not a 
finite extension of an abelian group.  But then neither is $G$ as $R$
is a quotient of $G$.

Conversely, assume $\chi (C,D)\geq 0$. If $(C, D)$ projectivizes to 
$(\PP^1, \bar D)$ and $\bar D$ has no more than two components, then
$G$ is abelian by
proposition~\ref{prop:elliptic-pi} and proposition~\ref{prop:star}. 
Otherwise, $(C,D)$ is uniformizable and, with the
notation as in section 3.3.1, the exact sequence \ref{seq} implies that
$G$ is an extension of $H$ by $R$.
Proposition~\ref{prop:A} now applies to give a pull back elliptic fibration
$\hat f:\hat X\ra \hat C$ without multiple fibers such that $\hat X$
is a finite unramified covering of $X$
and such that $\hat C=\PP^1$, $\CC^\ast$ or
an elliptic curve. We will consider each of these cases for $\tilde C$ 
separately. Note first that $G=\pi_1(X)$ (respectively $R$) is a finite
extension of $\hat G=\pi_1(\hat X)$ (respectively $\hat R$) and that
$\hat H=H$. Replacing $\hat C$ by a finite 
unramified covering of $\hat C$, we
may assume, thanks to lemma~\ref{lemma:IST}
and proposition~\ref{prop:kod}, that $H$ lies in the
center of $\hat G$ (that is, the conjugation action of $\hat G$ on $H$
is trivial).

In the case when $\hat C=\PP^1$, proposition~\ref{prop:elliptic-pi} implies
that $\hat G$ is a quotient of a free abelian group of rank two. 
Hence $\hat G$
is abelian of rank no greater than two. Since $X$ is K\"ahler if and 
only if $\hat X$ is, this rank is even if $X$ is  K\"ahler and odd
if not. 

In the case when $\hat C=\CC^\ast$, the triviality of the conjugation
action of $\hat R=\ZZ$ implies immediately that $\hat G$ is abelian,
of rank one greater than that of $H$.

In the case when $\hat C$ is an elliptic curve, 
proposition~\ref{prop:kod} implies that $\hat X$
must either be a holomorphically 
locally trivial fibration over $\hat C$, or
$H$ is finite cyclic. In the former case, proposition~\ref{prop:Ebundle}
tells us that $\hat G$ is a finite extension of a free abelian group of
rank four. In the latter case,
let $m/2$ be the order of $H$. Since $\hat R$
is abelian, the commutator of two elements $a$ and $b$ in $\hat G$ must
lie in $H$. Hence, $ab=bac$ for some $c\in H$. Since $c$  commutes with
both $a$ and $b$, we have $a^mb=ba^m$ and $a^mb^m=(ab)^m$.
This shows that $\hat G^m=\{a^m\ |\ a\in \hat G\}$ 
is an abelian subgroup of 
$\hat G$ intersecting $H$ at $1$. 
Hence, we can form the internal direct sum $G^m\oplus H$ in $G$
which we can easily identify with the inverse image of $\hat R^m$ in $G$,
where $\hat R^m$ is a normal subgroup of index  $m^2$ in $\hat R$. (We note as
an aside that $G^m$ is canonically isomorphic to $R^m$.)
It follows that $\hat G$ becomes
abelian if we replace $\hat X$ by a finite covering of itself and so 
our theorem is proved.
\BOX

\section{Other compact complex surfaces}

We deal with the remaining cases of compact complex surfaces in this section.
These are the case of zero Kodaira dimension and the case of zero algebraic
dimension. In fact, by
Kodaira's classification, all surfaces with Kodaira dimension zero are
elliptic fibrations except for  those bimeromorphic to compact complex 
2-dimensional tori and K3 surfaces, where the elliptic
ones form a dense codimension one family in their respective moduli space.
As we have already resolved the case of elliptic fibrations
in the previous section, we need
to consider only the tori and the K3 surface cases.
We first resolve the case of tori,
and indeed prove a much stronger result 
of independent interest, before considering the other cases.

\subsection{Compact complex tori}  \label{section:tame}

A 2-dimensional compact complex torus
is the quotient of $\CC^2$ by a lattice
$\Lambda$ of real rank 4.  
Let $X$ be such a surface, which we call a torus surface. 
Any compact surface $Y$ bimeromorphic to $X$
admits a dominating holomorphic map from the complement of 
finitely many points in $X$.
We show in this section that the complement of 
finitely many points in $X$
is dominable by $\CC^2$.  This will follow
immediately from proposition~\ref{prop:tame} below. Hence, $Y$ is
also dominable by $\CC^2$ as a result.

Following Rosay and Rudin \cite{RR}, we say that a discrete set
$\Lambda$ in $\CC^2$ is {\it tame} if there is a holomorphic automorphism,
$F$, of $\CC^2$ such that $F(\Lambda)$ is contained in a complex line.  Using
techniques of \cite{RR} or \cite{BF}, the complement of a tame set is
dominable by $\CC^2$, and in fact, there exists an injective
holomorphic map from $\CC^2$ to $\CC^2 \setminus \Lambda$. 

By a lattice, we mean a discrete $\ZZ$-module.  For the following
proposition, let $\Lambda$ be a lattice in $\CC^2$, let $q_1, \ldots, q_m
\in \CC^2$, and let $\Lambda_0= \cup_{j=1}^m \Lambda + q_j$, where
$\Lambda+q_j$ represents translation by $q_j$.  

\begin{prop}  \label{prop:tame}
The set $\Lambda_0$ is tame.  In particular, $\CC^2 \setminus
\Lambda_0$ is dominable by $\CC^2$ using an injective holomorphic map.
\end{prop}

This result will be strengthened considerably in
section~\ref{section:balls}.  Before proving this proposition, we need
a lemma.  

\begin{lemma}  \label{lemma:discrete}
There exists an invertible, complex linear transformation $A:\CC^2
\ra \CC^2$ such that ${\rm Im} \; \pi^1 A(\Lambda_0)$ is a discrete set
in $\RR$.  Moreover, we may assume that if $p,q \in A(\Lambda_0)$ with
$p \neq q$, then $|p-q| \geq 1$ and either ${\rm Im} \;\pi^1 p = {\rm
Im} \; \pi^1 q$ or $|{\rm Im} \; \pi^1 p - {\rm Im}\; \pi^1 q| \geq 1$.  
\end{lemma}

{\bf Proof:} Let $v_1, v_2, v_3, v_4$ be a $\ZZ$-basis for $\Lambda$,
and let $E$ be the span over $\RR$ of $v_1, v_2, v_3$.  Using the real
inner product, let $u_0 \neq 0$ be orthogonal to $E$.  Using the
complex inner product, let $u_1 \neq 0$ be orthogonal to $u_0$.
Then $u_1$ and $i u_1$ are both real orthogonal to $u_0$, so $\CC u_1
\subseteq E$.  Choose $A_1$ complex linear such that $A_1(u_0) = (1,0)$
and $A_1(u_1) = (0,1)$.  Then $\pi^1 A_1(E)$ is a one (real) dimensional
subspace of $\CC$, so by rotating in the first coordinate, we may
assume that $\pi^1 A_1(E)$ is the real line in $\CC$.  

Let $\mu_0 = {\rm Im}\; \pi^1 A_1(v_4)$, and $\mu_j = {\rm Im}\; \pi^1
A_1(p_j)$ for $j = 1, \ldots, m$.  Then for each $j = 1, \ldots, m$
and $k \in \ZZ$, we have
$$ \pi^1 A_1(E + kv_4 + p_j) \subseteq \RR + i(k\mu_0+\mu_j), $$
so that ${\rm Im}\;\pi^1 A_1(\Lambda_0)$ is discrete in $\RR$.  Applying an
appropriate dilation to $A_1$ gives $A$ as desired.  $\BOX$

Note that this lemma implies that given a finite
set of points in a complex 2-torus, there is an open set, $U$, containing
this finite set and a nonconstant image of $\CC$ avoiding $U$.
In particular, the complement of $U$ in this torus is not Kobayashi
hyperbolic.  As mentioned in the introduction, this result will be
strengthened in section~\ref{section:balls} to show that there is a
dominating map into the complement of such an open set $U$.

\bigskip
\noindent
{\bf Proof of proposition~\ref{prop:tame}:} Lemma~\ref{lemma:discrete}
implies that there is a complex line $L=\CC(z_0, w_0)$ with orthogonal
projection $\pi_L: \CC^2 \ra L$ and real numbers $\mu_0, \ldots,
\mu_m$ such that 
\begin{equation} \label{eqn:lines}
\pi_L(\Lambda_0) \subseteq \cup_{j=1}^m (\mu_0 \ZZ + \mu_j + i
\RR)(z_0, w_0).
\end{equation}
I.e., identifying $L$ with $\CC$ in the natural way, the image of
$\Lambda_0$ under $\pi_L$ is contained in a union of lines parallel to
the imaginary axis, and this union of lines intersects the real axis
in a discrete set.   

Making a linear change of coordinates, we may assume that $L =
\CC(0,1)$, in which case we may identify $\pi_L$ with projection to
the second coordinate, $\pi^2$.  Let $\pi^1$ denote projection to the
first coordinate, and let $E = \cup_{j=1}^m (\mu_0
\ZZ + \mu_j + i \RR)$.

\medskip
We next show that there is a continuous, positive function $f_0$ on
$E$ such that if $(z,w) \in \Lambda_0$ with $z\neq 0$, then $f_0(w)|z|
\geq 2|w|$.  First, define 
$$ r_1(w) = \left\{\begin{array}{ll}
\frac{|w|}{\min\{|z|: (z,w) \in \Lambda_0, z \neq 0\}}
        & {\rm if} \ \  w \in \pi^2(\Lambda_0); \\
0
        & {\rm if} \ \ w \in E \setminus \pi^2(\Lambda_0).
\end{array} \right. $$
Then $r_1(w) \geq 0$, and since $\Lambda_0$ is discrete, $r$ is
upper-semicontinuous.  

Let $r_2(w) = 2(r_1(w)+1)$ for $w \in E$.  Since $r_2$ is also 
upper-semicontinuous, it is bounded above on compacta, so a standard
construction gives a function $f_0$ which is continuous on $E$ with
$f_0(w) \geq r_2(w) >0$.  Then for $(z,w) \in \Lambda_0$ with $z \neq
0$, we have $f_0(w) |z| \geq 2 r_1(w)|z| \geq 2 |w|$ by definition of
$r_1$. 

\medskip
We next find a non-vanishing entire function $f$ so that $|f(w) z| \geq
|w|$ if $(z,w) \in \Lambda_0$ with $z \neq 0$.  
Since $f_0$ is positive on $E$, $\log f_0(w)$ is continuous and
real-valued on $E$, and $\log f_0(w)  \geq  \log 2 + \log (r_1(w) +
1)$.   By Arakelian's theorem (e.g. \cite{RR2}), there exists 
an entire $g(w)$ with $|\log f_0(w) - g(w)|< \log 2$ for $w \in E$.
Then $f(w) = \exp(g(w))$ is entire and non-vanishing, and if $(z,w) \in
\Lambda_0$ with $z \neq 0$, then 
$$ |f(w) z| = \exp({\rm Re} g(w)) |z|  \geq r_1(w)|z| \geq |w|. $$

\medskip
Finally, define $F(z,w) = (f(w)z, w)$.  Then $F$ is a biholomorphic
map of $\CC^2$ onto itself, and for $(z,w) \in \Lambda_0$ with $z \neq
0$, we have $|\pi^1 F(z,w)| \geq |\pi^2 F(z,w)|$.  Since
$F(\Lambda_0)$ is discrete, we see that $\pi^1 F(\Lambda_0)$ is
discrete.  Hence $F(\Lambda_0)$ is tame by \cite[theorem 3.9]{RR}.  By
definition of tame, $\Lambda_0$ is also tame, so as mentioned earlier,
$\CC^2 \setminus \Lambda_0$ is dominable by $\CC^2$.  $\BOX$

\begin{cor}\label{cor:tori} The complement of a finite set of points
in a two dimensional compact complex torus is dominable by $\CC^2$.
Hence any surface bimeromorphic to such a torus is dominable by $\CC^2$..
\end{cor}

We remark that not all tori are elliptic. The elliptic torus surfaces
form a $3$ dimensional family in the $4$ dimensional family of torus
surfaces and the generic torus contains no curves. All compact complex
tori are K\"ahler. Also a compact surface bimeromorphic to a torus can be
characterized by $\kappa=0$ and $q=2$.

\subsection{K3 surfaces} \label{subsec:k3}

A compact complex surface $X$ is called a K3 surface if its fundamental
group and canonical bundle are trivial. A useful fact in the compact complex
category, due to Siu (\cite{Siu}), is that all K3 surfaces are K\"ahler.
One can show that $H^2(X,\ZZ)$ is isometric to a fixed lattice $L$ of
rank $22$. If $\phi$ is such an isometry, then 
$(X,\phi)$ is called a marked K3 surface. The set of such surfaces is
parametrized by a $20$ dimensional non-Hausdorff manifold $\cal M$
\cite[VIII]{BPV} (The fact that $\cal M$ is smooth follows from S.T.
Yau's resolution of the Calabi conjecture in \cite{Yau}
(see e.g., \cite{T}) and
the fact that $\cal M$ is not Hausdorff is due
to Atiyah (\cite{At}).)

We first observe a few facts from the classical work of
Piatetsky-Shapiro and Shafarevich  in \cite{PS} 
(see also \cite{LP},\cite{Shi},\cite[VIII]{BPV}),
where they obtained a
global version of the Torelli theorem for K3 surfaces.
Given a marked K3 surface and a point 
$o\in {\cal M}$ corresponding to it, there is a smooth Hausdorff 
neighborhood ${\cal U}$ of $o$, a smooth complex manifold $Z$, and a
proper holomorphic map $Z \stackrel{p}{\rightarrow} {\cal U}$ whose
fibers are 
exactly the marked K3 surfaces parametrized by ${\cal U}$.
Within this local family, the subset of projective K3 surfaces
is parametrized by a topologically dense subset of
$\cal U$ which is a countable union of codimension one subvarieties. The
elliptic K3 surfaces (that is, K3 surfaces admitting an elliptic
fibration) also form a topologically dense codimension one family in
${\cal U}$.

The following proposition follows directly from theorem~\ref{thm:E2}
and the fact that the fundamental group of a K3 surface is trivial.

\begin{prop} \label{prop:elliptic-K3}
A compact complex surface bimeromorphic to an 
elliptic K3 surface is holomorphically dominable by $\CC^2$.
\end{prop}

The previous section on complex tori allows us to deal with another
class of K3 surfaces --- the Kummer surfaces, which form a 4
dimensional family in the 20 dimensional family of K3 surfaces. 
Such a surface $X$ is, by definition, obtained by taking the 
quotient of a torus surface $A$ (given as a complex
Lie group $\CC ^ {2}$/lattice) 
by the natural involution $g(x)=-x$, then blowing up the 16
orbifold singular points (resulting in 16 ($-2$) curves). Alternatively,
one can describe $X$ as a $\ZZ_2$ quotient of $\hat A$, where $\hat A$ is the
blowing up of $A$ at the 16 points of 
order $2$ and where the quotient map is branched
along the exceptional (-1)-curves of the blowing up.
Since the inverse image of any finite set of points in $X$ is finite
in $\hat A$ and hence also finite in $A$, any surface bimeromorphic
to a Kummer surface is dominable by $\CC^2$ according to 
corollary~\ref{cor:tori}.

\begin{prop}\label{prop:Kummer-K3}
A compact surface bimeromorphic to a Kummer surface is dominable by $\CC^2$.
\end{prop}

Before we leave the subject of K3 surfaces, it is worth mentioning
that projective K3 surfaces are dominable by $\DD\times\CC$ by the work
of \cite{GG} and \cite{MM}. Clearly, elliptic K3 surfaces 
and Kummer surfaces are so dominable
as well. Such a surface cannot be measure hyperbolic as defined by
Kobayashi (\cite{Kob}).
However, it is still an unsolved problem
whether all K3 surfaces are so dominable. The only other compact complex
surfaces for which this problem remains open are the
non-elliptic and non-Hopf surfaces of
class VII$_0$ outside the Inoue-Hirzebruch construction.

\subsection{Other compact surfaces and our two main theorems}

Besides those bimeromorphic to
K3 and torus surfaces, the remaining compact complex surfaces
with zero Kodaira dimension are all elliptic, 
and are all dominable by $\CC^2$.
Such a surface must be bimeromorphic to either a Kodaira surface (defined
and characterized in section~\ref{sec:pi_1}), a hyperelliptic surface
(which is a finite free quotient of a product of elliptic curves,
and hence projective),
or an Enriques surface (which 
is a surface admitting an unramified double covering by
an elliptic K3 surface). Except for the first among these three types,
the fundamental group is always a finite extension of an abelian group.

\smallskip

Finally, the only remaining compact complex surfaces
are those with algebraic dimension 0 and $\kappa =-\infty$. 
This category includes the non-elliptic Hopf surfaces, which are
dominable by $\CC^2$ by construction (see \cite{Ko4}).  This category
also includes the Inoue surfaces, which must be excluded from our main
theorems since their universal cover is $\DD\times\CC$, hence are
not dominable by $\CC^2$, while any nonconstant image of $\CC$ must be
Zariski dense (see proposition 19.1 in \cite[V]{BPV}).
However, it is of interest to note that the Zariski dense 
holomorphic images of $\CC$ are constrained
by higher order equations on an Inoue surface so that if we relax
property C in this sense, we can in fact 
include Inoue surfaces in the next theorem.
Unfortunately, aside from the Hopf surfaces 
and the Inoue surfaces, the detailed structure of
surfaces of this type is not yet clear even though we know 
the existence of projective affine structures for a special
subclass of these surfaces.

\medskip

We now summarize our investigation in the compact category
by giving the following extensions of our main theorems stated in the
introduction:

\begin{theorem} \label{thm:1'}
Let $X$ be a compact complex surface of Kodaira dimension less than $2$.
Assume that either $\kappa(X)\neq -\infty$ or $a(X)\neq 0$. In the case
that $X$ is bimeromorphic to a K3 surface that is not Kummer, 
assume further that $X$ is elliptic.
Then $X$ is dominable by $\CC^2$ if and only if it does not satisfy
property C.  Equivalently, there is a 
dominating holomorphic map $F:\CC^2 \ra X$
if and only if there is a holomorphic image of
$\CC$ in $X$ which is Zariski dense.
\end{theorem}

\begin{theorem} \label{thm:2'}
Let $X$ be a compact complex surface not bimeromorphic to a Kodaira surface.
Assume that either $\kappa(X)\neq -\infty$ or $a(X)\neq 0$. In the case
that $X$ is bimeromorphic to a K3 surface that is not Kummer, 
assume further that $X$ is elliptic.
Then $X$ is dominable by $\CC^2$ if and only if it has Kodaira dimension 
less than two and its fundamental group is a finite extension of an
abelian group (of rank $4$ or less). 
\end{theorem}


\section{Non-compact algebraic surfaces}
\label{sec:noncompact}

We begin with a  key example which motivated
the general algebraic setting. This is the example of the complement
of a smooth cubic curve in $\PP^2$, which we will show to be dominable
by $\CC^2$.

\subsection{Complement of a cubic in $\PP^2$}  \label{subsec:cubic}

Let $C$ be a smooth cubic curve in $\PP ^ {2}$ and let $X = \PP^2
\setminus C$.  Then its logarithmic canonical bundle $K_{\PP^2}(D)$
is the trivial line bundle as $\deg K_{\PP^2}=-3$.
Hence, $\bar{\kappa} (X) = 0$ and $X$ is a logarithmic
K3 surface; that is, a non-compact 2-dimensional Calabi-Yau manifold.




\begin{prop}\label{prop:cu}
The surface $X = \PP^2 \setminus C$ is dominable by $\CC^2$.
\end{prop}

\noindent
{\bf Proof:} A tangent line to $C$ at a non-inflection point meets $C$
at one other point.  This gives rise to a holomorphic $\PP ^{1}$ bundle
with two holomorphic sections. To see that this is actually a bundle
(i.e. locally trivial), identify it with the projectivization of the
tautological vector bundle of rank two over the dual curve of $C$ with
the obvious isomorphism. We may pull back this $\PP ^{1}$ 
bundle and the sections
to the universal cover $\CC$ of $C$, with two sections $s_{\infty}$ and
$s$.  Hence one may regard the complement of $s_{\infty} ( \CC )$ of
this bundle as a trivial line bundle on $\CC$ with a meromorphic
section $s$ (with poles coming from points of inflection of the cubic).

Hence, it suffices to construct a holomorphic map from $\CC ^ {2}$
onto the complement of the graph of a meromorphic function $s$ to
give a dominating map to $X$.
Note that each vertical slice 
of the complement of the graph is $\CC ^ {\ast}$ except at a pole of $s$,
where the vertical slice is $\CC$. 

To construct such a map, first define
\begin{eqnarray}\label{eqn:graph}
\psi (t,w) & = & \frac{\exp (tw) -1}{t} \\
& = & w + \frac{tw ^ {2}}{2!} + \frac{t ^ {2} w ^ {3}}{3!} + \cdots
\end{eqnarray}
which is entire on $\CC ^ {2}$. Note that $(t,\psi (t,w))$ 
is a fiberwise selfmap of $\CC^2$ which misses precisely
the graph of  $-1/t$, a function with a simple pole
at the origin.  

Since $\CC$ is Stein, there exists an entire function $g$ such
that $\frac{1}{g}$ has the same  principle parts as $s$.  This is 
because we may write $s=f/f_1$ where $f$ and $f_1$ are entire
with no common zeros. So $\log f$ is well defined in a neighborhood of
each zero of $f_1$. By Mittag-Leffler and Weierstrass, we can find an
entire function $g_1$ with the same Taylor expansion as $\log f$ 
to the order of vanishing of $f_1$ at each zero of $f_1$. Then 
$g=f_1/\exp g_1$ is our desired function. In
particular, $g$ vanishes precisely when $s$ has a pole.
Then $h = s - \frac{1}{g}$ is entire, so 
\begin{eqnarray*}
\phi (z,w) & = & h(z) - \psi (g(z), w) \\
           & = & s(z) - \frac{\exp(w g(z))}{g(z)}
\end{eqnarray*}
is entire on $\CC ^ {2}$.  For fixed $z$ with $g(z) \neq 0$, we see
from the second equality that $\phi(z,w)$ can attain
any value in $\CC \setminus \{s(z)\}$  by varying $w$.  If $g(z) = 0$, then
$\phi(z, w) = h(z) - w$, which can attain any
value in $\CC$ by varying $w$.  

Hence, the map $\Phi: \CC^2 \rightarrow \CC^2 \setminus
\mbox{graph}(s)$ given by 
$$\Phi(z,w) = (z, \phi(z,w))$$
is holomorphic and onto.  Composing this map with the map into the
$\PP^1$ bundle over $C$, we obtain a dominating map into the
complement of the cubic $C$.  $\BOX$

Note that an important step here is the construction of an entire function
$h$ whose graph does not intersect the graph of $s$. This is certainly
analogous to the situation of elliptic fibrations.

\medskip \noindent

{\bf Remark:} The complement of a smooth cubic
does not admit any algebraic map to $\PP^1$ whose generic fiber 
contains $\CC^\ast$. This is
the only example among complements of normal crossing 
divisors in $\PP^2$ with this property. In fact, this is the 
only meaningful affine example with this property that is dominable by
$\CC^2$ (see \cite[p. 189]{M}).  Since this
is a logarithmic K3 surface, this phenomenon is suggestive
of the situation for a generic compact K3 surface.

\medskip \noindent 
We isolate the following useful theorem from the above proof.

\smallskip

\begin{theorem} \label{cubic}
Let $s$ be a meromorphic function on $\CC$. Then the
complement of the graph of $s$ admits a 
dominating fiber-preserving holomorphic map 
from $\CC^2$. 
\end{theorem}

\smallskip

\subsubsection{Complements of normal crossing divisors in $\PP^2$}

Let $X$ be the complement of a normal crossing divisor $D$ in $\PP^2$.
If $\deg D>3$, then $\bar \kappa (X)=2$ and hence $X$ is not dominable by
$\CC^2$. If $\deg D=3$, then $D$ consists of at most three components
and it is easy to check that $X$ is dominable by $\CC^2$ as follows. 
If $D$ has only one component, then it is either a smooth cubic or a
cubic with one node. In the first case, the result
follows from proposition~\ref{prop:cu}. In the second case, blowing
up that node gives us a $\PP^1$ bundle over $\PP^1$ with two sections,
one corresponding to the exceptional curve of the blow-up. These
two sections intersect precisely at the two fibers of the bundle
corresponding to the two tangent directions of the cubic at the node.
Hence, removing these two fibers gives us a surface biholomorphic to
$\CC^\ast\times\CC^\ast$, which is dominable by $\CC^2$.
If $D$ has two components, then it
consists of a line and a conic (that is, a smooth curve of degree two)
intersecting at two points.  Blowing up one of the  
points of their intersection (corresponding to projecting from this point
of intersection) gives us a $\PP^1$ bundle over $\CC$ with two
sections complemented, one of which is the exceptional
curve of the blow-up.  If we think of one section as $\infty$, then the
other section can be regarded as a meromorphic function on $\CC$ and so
theorem~\ref{cubic} applies to give a dominating map from $\CC^2$ to $X$.
An easier way is to delete the fiber containing the only point of intersection
of these two sections. The resulting $X$ is biholomorphic to 
$\CC^\ast\times\CC^\ast$ and hence dominable by $\CC^2$.
If $D$ has three components, then each
must be a line and $X$ is  $\CC^\ast\times\CC^\ast$, which is dominable
by $\CC^2$. 

{}From the above argument, we see also that if $\deg D<3$, then $X$ is
dominable by $\CC^2$. In summary, we have:

\begin{theorem} Let $D$ be a normal crossing divisor in $\PP^2$.
Then $\PP^2\setminus D$ is dominable by $\CC^2$ if and only if
$\deg D\leq 3$.
\end{theorem}

We remark that this theorem is no longer true if $D$ is not normal 
crossing. The unique counterexample in one direction is when $D$ consists of
three lines intersecting at only one point, which is not dominable by $\CC^2$. 
Another counterexample, but in
the opposite direction, is given by the complement of the union of a conic
and two lines intersecting at a point of the conic (which we discussed in the
two component case of $\deg D=3$ above).

\subsection{The general quasi-projective case}

Let $X$ be an algebraic surface over $\CC$. Then 
$X=\bar{X}\setminus D$ 
where $\bar{X}$ is projective and $D$ is a normal
crossing divisor in $\bar{X}$. This is the notation set 
forth in section 2 and we will assume this setup throughout this section.
Kawamata (\cite{K1},\cite{K2},\cite{K3}) has considered the structure of $X$
and obtained a classification theory analogous to that in the projective
case. Much of this is explained in 
some detail in Miyanishi (\cite{M}). We
will use their results directly to tackle our problem in this section.

If there is a surjective morphism  $f:X\ra C$
whose generic fiber is connected,  then we say that $X$ is
fibered over $C$.  (We remind the reader that morphisms are algebraic
holomorphic maps.) More generally, if $f$ is required to be only
holomorphic rather than a morphism, then we say that $X$ is 
holomorphically fibered over $C$. For example, the complement of
the graph of a meromorphic function is holomorphically fibered over
$C$ with generic fiber $\CC^\ast$.
As before, we let $X_s=f^{-1}(s)$ be the fiber over $s$.
We first quote the subadditivity property of 
(log-)Kodaira dimension due to Kawamata (\cite{K1}):

\begin{prop}\label{prop:sub}
If $X$ is fibered over a curve $C$, then
$$\bar\kappa(X)\geq \bar\kappa(C)+\bar\kappa(X_s)$$
for $s$ outside a finite set of points in $C$; that is, for 
the generic fiber $X_s$.
\end{prop}

$\mbox{}$From the definitions, 
a curve of positive genus with punctures has positive
Kodaira dimension. An elliptic curve has Kodaira dimension zero.
A punctured $\PP^1$ has $\kappa=-\infty, 0$ or $1$
according to the number of punctures being $1,2$ or greater than $2$,
respectively.

Given a dominating morphism $f$ between
algebraic varieties, 
it is clear  that $f^*$ is injective on the level of logarithmic
forms (see \cite{Ii}). 
Since tensor powers of top dimensional logarithmic forms define
the Kodaira dimension, we see that if $f$ is equidimensional, then it
must decrease Kodaira dimension.

If $\bar q(X)>0$, then there is a morphism from $X$ to a 
semi-abelian variety
(a commutative algebraic Lie group that is an extension of a 
compact torus by $(\CC^\ast)^k$ for some $k$)
of dimension $\bar q(X)$, called the quasi-Albanese
map and constructed by Iitaka in \cite{Ii1}. One has the simple
formula relating $\bar q(X)$ to the first Betti numbers of $X$ and $\bar X$:
$$\bar q (X) - q (\bar X) = b_1 (X) -b_1(\bar X).$$

Note that $\CC$ does not support any logarithmic form by this formula.

\subsubsection{Surfaces fibered by open subsets of $\PP^1$}

Let $X$ be fibered over a curve C by a map $f$ whose generic fiber is
$\PP^1$ (possibly) with punctures. Then,  by a finite number of
contractions of $(-1)$-curves that remain on the fiber, 
the compactification $\bar X$ of $X$ admits a birational
morphism $g$ to a ruled surface $\bar Y$ over a projective curve $\bar C$, the
compactification of $C$, and $g$ is a composition of blowing ups.
Hence $Y=\bar{Y}|_C$ is a $\PP^1$ bundle over $C$, whose bundle map will
again be denoted by $g$. We may write $f=h\circ g$, where
\begin{equation}\label{gh}
\CC \stackrel{r}{\ra} X \stackrel{g}{\ra} Y \stackrel{h}{\ra} C.
\end{equation}

If every holomorphic image of $\CC$ in $X$ is constant in $C$ (when 
composed with $f$), then $X$ satisfies property C.
Otherwise, there exists a holomorphic map $r:\CC\ra X$ such that 
$f\circ r$ is not constant. By taking the fiber product with $f\circ r$,
we can pull back the factorization picture \ref{gh} to one over $\CC$
$$\CC \stackrel{\tilde{r}}{\ra} \tilde X \stackrel{\tilde g}{\ra}
\tilde Y\stackrel{\tilde h}{\ra} \CC,$$
where  $\tilde f=\tilde h\circ \tilde g$ is surjective with a holomorphic
section $\tilde r$. Here, $\tilde X$ may be singular, but we will 
regard it only as an auxiliary space.

\medskip

We will first deal with the case where the general fiber has at most
one puncture; that is, $X_s=\PP^1$ or $\CC$ for $s$ in an open subset 
of $C$.
We can then regard $\tilde Y$ as a trivial $\PP^1$ bundle with a section
$D_\infty$ to which the puncture (if one exists) on the ``generic'' fiber of
$\tilde f$ is mapped.  Note that $\tilde Y\setminus D_\infty=\CC^2$
with coordinates $(z,w)$, and so we may regard a section of $\tilde h$
as a meromorphic function on $\CC$. In particular, $\tilde g \circ
\tilde r$ is a meromorphic section of $\tilde h$. 

Since $\bar X$ is obtained from $\bar Y$ by a finite number
of blow ups, we can identify points on $X$ as infinitely near points
on $Y$ of order 0 or more as in \cite[p. 392]{Ha}. Note that the 
set of fibers in $Y$ which contain infinitely near points of
order 1 or more is finite (since the set of such fibers in $\bar
Y$ is finite). This finite set of fibers in $Y$ pulls back to a
discrete set of fibers in $\tilde Y$.  In $Y$, such a higher order
infinitely near point corresponds to a point in $X$ obtained by
finitely many blow-ups, hence to the specification of a finite jet at
the point in $Y$.  Under pull-back, this corresponds to a finite jet
in $\tilde Y$.  Additionally, there is a finite set of fibers in $Y$
which may have more than one puncture, and these fibers all pull back
to a discrete set of fibers in $\tilde Y$.  Together, these two types
of fibers will be called exceptional fibers.  

In order to produce a dominating map into $X$, it suffices to produce
a fiberwise dominating map $F(z,w) = (z, H(z,w))$ into $\tilde Y$
which respects these exceptional fibers in the following sense.  If
$\tilde{Y}_s$ is an exceptional fiber, then $F(s,w)$ is a single
point independent of $w$.  Moreover, if $\tilde{Y}_s$ is a fiber having
more than one puncture, then the image of the map $F$ should avoid all
such punctures.  If $\tilde{Y}_s$ is a fiber having an infinitely near
point, then $F(s,w)$ should equal $\tilde g \circ \tilde r(s)$.
Additionally, if $\tilde g \circ \tilde r$ passes through this
infinitely near point, and $\phi$ is holomorphic in a neighborhood of
$s$, then the local curve $z \mapsto (z, H(z,\phi(z)))$ should agree
with the jet given by the infinitely near point on $\tilde{Y}_s$.

Fortunately, the section $\tilde g \circ \tilde r$ has the correct jet
whenever it intersects one of these exceptional fibers, so we can use
this section to obtain such a map.  Let $q(z) = \tilde g \circ \tilde
r(z)$, which is meromorphic.  We will define $H(z,w) = p(z) w + q(z)$ for
some entire $p(z)$.  For each exceptional fiber $\tilde{Y}_s$, there
is an integer $n_s \geq 1$ such that if $p$ vanishes to order $n_s$ at
$s$, then $F$ defined with this $H$ respects the exceptional fiber as
indicated above.  By Weierstrass' theorem, there exists $p$ entire
vanishing exactly to order $n_s$ at each $s$.  Then $F(z,w) = (z,
H(z,w))$ gives a dominating map from $\CC^2$ into $\tilde Y$
respecting the exceptional fibers, and this map pushes forward to $Y$,
then lifts to give a dominating map into $X$, as desired.

\smallskip

We now deal with the case where the generic fiber of $f$ is $\CC^*$. 
In this case, $\tilde Y$ can be identified with a $\PP^1$
bundle with a double section $D_Y$, to which the punctures on the
``generic'' fibers of $\tilde f$ maps to. Now, either $D_Y$ consists
of two components, both of which are smooth sections of $\tilde h$, or
$D_Y$ consists of one component. In either case, outside of a discrete
set of fibers, $D_Y$ can be written locally as the union of two 
mermorphic sections.  Moreover, we define the set of exceptional
fibers exactly as in the previous case. 

First, using a fiber-preserving biholomorphic map of $\CC \times
\PP^1$ to itself, we may move $\tilde g \circ \tilde r$ to become the
$\infty$-section. Then the requirement of
agreeing with the jet of $\tilde g \circ \tilde r$ at a point $s$ is
equivalent to having a pole of some given order at $s$ in the new
coordinate system.  Next, let $E_1$ be the
points in $\CC$ at which $D_Y$ intersects this new
infinity section.  Near a point $s \in E_1$, $D_Y$ can be written as $w
= h(z) \pm \sqrt{g(z)} = u^\pm(z)$ for some meromorphic $g$ and $h$.
Hence there exists $n_s > 0$ such that $u^\pm(z) (z-s)^{n_s}$
converges to 0 as $z$ tends to $s$.  We may assume also that if $s \in
E_1$ and $s$ is the base point of an exceptional fiber, then the $n_s$
obtained here is larger than the $n_s$ obtained above for this exceptional
fiber. 

Let $E$ be the union of $E_1$ and the set of base points corresponding
to exceptional fibers.  Let $p$ be entire with a zero of order $n_s$
at each $s \in 
E$ and no other zeros, and let $\Phi(z,w) = (z, p(z)w)$.  Then
$\Phi(D_Y)$ is a double section in $\CC \times \PP^1$, and a dominating
map from $\CC^2$ to $\CC^2 \setminus \Phi(D_Y)$ followed by
$\Phi^{-1}$ gives a dominating map to the complement of $D_Y$ which
respects the exceptional fibers.  

Hence it suffices to construct a dominating map into the complement of
$\Phi(D_Y)$.  Note that $\Phi(D_Y)$ can be written as $w = v^\pm(z) =
p(z) u^\pm(z)$, where $v^\pm$ are holomorphic except possibly for
square root singularities at branch points.  

For complex numbers $u$ and $v$, define a Mobius transformation $N_{u,v}(w)
= (uw-v)/(w-1)$, which takes 0 to $v$ and $\infty$ to $u$, and define
$G_{u,v}(w) = \exp(w(u-v))$.  Note that $N_{u,v}(w) = N_{v,u}(1/w)$
and that $G_{u,v}(w) = 1/G_{v,u}(w)$.  Hence $H_0(u,v,w) =
N_{u,v}(G_{u,v}(w))$ satisfies $H_0(u,v,w) = H_0(v,u,w)$.  Since 
symmetric functions of $v^+$ and $v^-$ are holomorphic, we see that
$H(z,w) = H_0(v^+(z), v^-(z), w)$ is well-defined and holomorphic from
$\CC^2$ to $\CC \times \PP^1$.  Moreover, for fixed $s$ such that
$v^\pm(s)$ are distinct, $H(s,\cdot)$ is nonconstant from $\CC$
to $\PP^1 \setminus \{v^\pm(s)\}$.  If $v^\pm(s)$ are equal, then
assuming without loss that $s=0$, we have $v^\pm(z) = h(z) \pm \sqrt{g(z)}$
for some holomorphic $g(z) = z^m g_1(z)$ with $g_1(0) \neq 0$, $m \geq
1$.  Then $v^+ - v^- = 2 \sqrt{g}$, so multiplying the numerator and
denominator of $H$ by $\exp(-w(v^+-v^-)/2)$ and using the Taylor
expansion of $\exp$ gives
\begin{align*}
H &= \frac{(h+\sqrt{g})(1+w\sqrt{g}) - (h-\sqrt{g})(1-w\sqrt{g}) + O(|z|^m)}
 	{(1+w\sqrt{g}) - (1-\sqrt{g}) + O(|z|^m)}\\[2mm]
& = \frac{2hw\sqrt{g} + 2\sqrt{g} + O(|z|^m)}{2w\sqrt{g} + O(|z|^m)}.
\end{align*}
As $z \rightarrow 0$, this last expression tends to $h(0) + 1/w$, and
hence $H(0, \cdot)$ maps $\CC$ onto $\PP^1 \setminus \{v\pm(z)\}$.  

Thus $H$ is a dominating map from $\CC^2$ to the complement of $\Phi(D_Y)$,
hence as noted before, $\Phi^{-1} \circ H$ is a dominating map from
$\CC^2$ to the complement of $D_Y$ which respects the exceptional
fibers.  As before, this map pushes forward to $Y$ and lifts to give a
dominating map into $X$, as desired.  

\smallskip

We can now summarize with the following theorem.

\begin{theorem}\label{thm:punc}
Assume that $X$ is fibered over a curve $C$ and that
the generic fiber is $\PP^1$ with at most two punctures. Then
$X$ is dominable by $\CC^2$ if and only if there is a Zariski
dense image of $\CC$ in $X$.
\end{theorem}

The arguments given in this paper are not sufficient to resolve the
question of dominability for open fibered surfaces.  As an example, we
have the following question.
\begin{question}  Let $X$ be the complement of a double section in a
conic bundle over $\CC$, $\CC^*$, or an elliptic curve.  
Is $X$ dominable by $\CC^2?$
\end{question}
We will consider this and related questions in a forthcoming paper.  

\medskip

\subsubsection{The $\bar \kappa= -\infty$ case}

Let $\bar\kappa(X)=-\infty$. Then $\kappa(\bar X)=-\infty$ as well. 
Hence $\bar X$ is either rational or birationally ruled over a curve of 
non-negative genus. In the latter case, proposition~\ref{prop:sub}
says that $X$ is fibered over a curve $C$ with $\kappa (C)\geq 0$ where
the generic fiber is $\PP^1$ with at most one puncture.
Hence theorem~\ref{thm:punc} applies in this case to give us the
equivalence of dominability by $\CC^2$ and the failure of property
C. Note that property C holds in the case $\kappa (C)> 0$
(which include the case $\bar q(X)\geq 2$), 
corresponding to $C$ being hyperbolic.

In the remaining case when $\bar X$ is rational, we can again divide
into two cases according to whether $\bar q(X)$ is zero or not.
In the latter case, we again have a fibering of $X$ over a curve $C$
via the quasi-Albanese map
with the generic fiber having at most one puncture by 
proposition~\ref{prop:sub}, as before. This is because there are
no logarithmic $2$-forms on $X$ since $\bar\kappa(X)=-\infty$. By
the same token, every logarithmic $1$-form on $X$ is the pull
back of a logarithmic form on $C$ (One can also see this from the fact
that $\PP^1$ with at most one puncture has no logarithmic forms so that
any logarithmic form on $X$ becomes trivial when restricted to the 
generic fiber. Hence, $\bar q (X)=\bar q(C)$.)
So, $C$ must be $\PP^1$ with at least 
two punctures. If it has more than two punctures, 
corresponding to $\bar q(X)\geq 2$, then $C$ is hyperbolic. 
So we have degeneracy of holomorphic maps
from $\CC$ in this case. Otherwise, theorem~\ref{thm:punc} applies.

We are left with the case where $\bar q(X)=0$ where 
proposition~\ref{prop:sub} no longer applies,
but where much of the
analysis has been done in \cite{M}. We now quote theorem $(1')$ of \cite{M},
(which follows from theorem I.3.11 of \cite{M})
\begin{theorem}\label{thm:conn}
With $X$ and $D$ as before, assume that $D$ is connected. Then 
$\bar \kappa(X)=-\infty$ if and only if $X$ fibers over a curve with
generic fiber being $\PP^1$ or $\CC$.
\end{theorem}

Except in the case where $X=\PP^2$,
there is, of course,  some fibering of $X$ to 
a curve (as is clear from, for example, (1) of
the classification list given in section 2)
and every such fibering must be to a curve $C$ that is either
$\PP^1$ or $\CC$. In these fibered cases, we would like to show
that the generic fiber is $\PP^1$ with at most two punctures so that
theorem~\ref{thm:punc} can be applied to
show that $X$ is dominable by $\CC^2$.
However, it remains an open question whether or not the generic fiber
has this form, and although this question should be resolved by some
case checking, this lack prevents us from giving a complete
classification in the case $\bar \kappa (X)=-\infty$
and $\bar q(X)=0$.

	We can now summarize this section as follows.

\begin{theorem}\label{thm:infty}
Let $X$ be the complement of a normal crossing divisor $D$ in a projective
surface. Assume $\bar \kappa (X)=-\infty$. If $\bar q(X)\geq 2$, then $X$
satisfies property C and hence is not dominable by $\CC^2$.
If $\bar q(X)=1$ or if $\bar q(X)=0$ and
$D$ is connected, then $X$ is dominable by $\CC^2$ if and only if
there exists a holomorphic map of $\CC$ to $X$ whose
image is Zariski dense. 
\end{theorem}

\subsubsection{The $\bar \kappa =1$ case}

Here, we can directly apply the basic Iitaka fibration theorem,
theorem 11.8 in \cite{Ii} (see also \cite{Ue}):

\begin{theorem} \label{thm:Ii}
Assume $\bar \kappa (X)\geq 0$. Then $X$ is properly
birational to a variety $X^*$ which is fibered over a variety of
dimension $\bar \kappa (X)$ and whose generic fiber has Kodaira 
dimension zero.
\end{theorem}

This theorem holds for $X$ of any dimension. But for our situation 
at hand, it says that $X$ is properly birational to a surface 
$X^*$ which is fibered over a curve with generic fiber
that is either an elliptic curve, or $\PP^1$ with two punctures.
Now, we have already shown that for such a fibered variety, 
dominability is unchanged for any variety properly birational to it.
The latter case is already resolved by theorem~\ref{thm:punc}.
The former case can also be resolved to give the same
conclusion by the same analysis as that of
theorem~\ref{thm:punc} with the help
of the Jacobian fibration as in section 3.
Thus, combining with theorem~\ref{thm:punc}, we have:

\begin{theorem}\label{thm:one}
Assume $X$ is fibered over a curve with generic fiber that is either
an elliptic curve or $\PP^1$ with at most two punctures. This is
the case, for example, when $\bar\kappa (X)=1$. 
Then $X$ is dominable by $\CC^2$ if and only if
there exists a holomorphic map of $\CC$ to $X$ whose
image is Zariski dense. 
\end{theorem}

\subsubsection{The $\bar \kappa =0$ case}

It remains to look at the case where 
$\bar \kappa(X)=0$. If $\bar q(X)\geq 2$, then a well known theorem
of Kawamata (\cite{K4}) says that $X$ has a birational morphism to
a semi-abelian surface. Hence, $X$ is dominable by $\CC^2$. 
If $\bar q(X)=1$, then $X$ is
fibered over a curve and the generic fiber is an elliptic curve or
is $\PP^1$ with at most two punctures by proposition~\ref{prop:sub}. 
Hence theorem~\ref{thm:one}
applies in this case. 
When $\bar q(X)=0$, our problem remains with 
some K3 surfaces as explained in section~\ref{subsec:k3}.

Finally, if $X$ is affine rational and 
$D$ has a component that is not a rational curve, then Lemma II.5.5
of \cite{M} says that either $X$ is fibered over a curve with
generic fiber $\PP^1$ with at most two punctures or $X$ is the
complement of a smooth cubic in $\PP^2$. The former is handled by 
theorem~\ref{thm:punc} while the latter is dominated by $\CC^2$ 
as shown in section~\ref{subsec:cubic}.
This resolves the case of the complement of a reduced divisor $C$ in
$\PP^2$ unless $C$ is a rational curve, which one can resolve
as well when $C$ has either only one singular point or is of low degree
(and it is easy to check all the cases for degree less than $4$).
This is a good exercise for the case when $C$ is a rational curve of
high degree, which we will not attempt here.
Note that, if $C$ is normal crossing with dominable complement, then
$C$ is again a smooth cubic in $\PP^2$, being the unique non-rational
component.

\begin{theorem}\label{thm:zero}
Assume $\bar \kappa (X)=0$. If  $\bar q(X)$ is positive or if $X$
is affine and $D$ has a component that is not a rational curve, then
$X$ is dominable by $\CC^2$ if and only if
there exists a holomorphic map of $\CC$ to $X$ whose
image is Zariski dense. 
\end{theorem}

\section{Compact complex surface minus small balls}\label{section:balls}

For the compact complex surfaces which we showed to be dominable by
$\CC^2$, a surprisingly stronger result can be achieved, thanks to 
the theory of Fatou-Bieberbach domains. We can show that these surfaces
remain dominable
after removing any finite number of sufficiently small open balls.
In this section we show how this can be done in the most difficult case,
the case of a two dimensional compact complex torus. 
We show that given any finite set of points in a
torus $T$, it is possible to find some open set, $U$, containing this
finite set, and a holomorphic map $F:\CC^2 \ra T \setminus U$
with non-vanishing Jacobian determinant.  In fact, $F$ lifts to an
injective holomorphic map from $\CC^2$ to $\CC^2$.  For the statement
of the following theorem, we focus only on this lifted map.  For
notation, $\Delta^2(p; r)$ is the bidisk with center $p$ and radii 
$r$ in both coordinate directions and $\pi^j$ represents projection to
the $j$th coordinate axis.  

\begin{theorem}  \label{thm:main}
Let $\Lambda \subseteq \CC^2$ be a discrete lattice, let $p_1, \ldots, p_m \in
\CC^2$, let $\Lambda_0 = \cup_{j=1}^m \Lambda + p_j$, and for $r>0$,
let $\Lambda_{0,r} = \cup_{p \in \Lambda_0} \Delta^2(p; r)$.  For some $r
>0$, there exists an injective holomorphic map $F: \CC^2 \ra
\CC^2 \setminus \Lambda_{0,r}$.  
\end{theorem}

In fact, the proof will show that any discrete set contained in
$\Lambda_{0,r}$ is a tame set in the sense of
section~\ref{section:tame}.  As an immediate corollary, we obtain the
following result, as mentioned in the introduction.  An
$n$-dimensional version of this result is found in \cite{Bu}.

\begin{cor}
Let $T$ be a complex 2-torus and let $E \subset T$ be finite.  Then
there exists an open set $U$ containing $E$ and a dominating map from
$\CC^2$ into the complement of $U$.  
\end{cor}

For the remainder of this section, $\Lambda$, $\Lambda_0$ and
$\Lambda_{0,r}$ will be as in the statement of this theorem.

\subsection{Preparatory lemmas}

In this subsection we state some necessary lemmas.  The proofs are
straightforward and perhaps even standard, but they are provided for
completeness.  

\paragraph{Notation:} For $\eps>0$, let $S_\eps = \{x+iy: x \in \RR,
-\eps < y < \eps\}$.  

\begin{lemma}  \label{lemma:g-int}
Let $C>0$, let $f:\RR \ra [0,C]$ be measurable, and let $\eps \in
(0,1)$.  Then there exists a function $g$ holomorphic on $S_\eps$ such
that if $\delta>0$ and $z_0 = x_0 + iy_0 \in S_\eps$ with $f(x) = c_0$
for $x_0 - \delta < x < x_0 + \delta$, then 
$$ |g(z_0) - f(x_0)| \leq \frac{2C\eps}{\pi \delta}. $$
Moreover, ${\rm Re} \; g(z) \geq 0$ for all $z \in S_\eps$.  
\end{lemma}

{\bf Proof:} For $n \in \ZZ$, let 
\begin{align*} 
g_n(z) &= \frac{1}{2\pi i} \int_{-n}^n \left( \frac{f(x)}{x-i\eps -z} -
\frac{f(x)}{x+i\eps -z} \right) dx \\
&= \frac{1}{\pi} \int_{-n}^n f(x) \left(\frac{\eps}{(x-z)^2 +
\eps^2}\right) dx. 
\end{align*}
I.e., $g_n$ is obtained via the Cauchy integral using the function $f$
on the two boundary components of $S_\eps$ and truncating at $x=\pm
n$.  By \cite[Thm 10.7]{rudin}, each $g_n$ is holomorphic in $S_\eps$.
Moreover, for $z_0 = x_0 + iy_0 \in S_\eps$, we have $|y_0| < \eps$, so
\begin{equation}  \label{eqn:bound}
|(x-z)^2+\eps^2| \geq {\rm Re}\;(x-(x_0+iy_0))^2 + \eps^2 \geq
(x-x_0)^2.
\end{equation}

Using this last inequality and the boundedness of $f$, it follows
immediately that $g_n$ converges uniformly on compact subsets of
$S_\eps$ to the holomorphic function  
\begin{equation}  \label{eqn:g-int}
g(z) = \frac{1}{\pi} \int_{-\infty}^\infty f(x)
\left(\frac{\eps}{(x-z)^2 + \eps^2}\right) dx. 
\end{equation}
A simple contour integration shows that if $f$ is replaced by the
constant $c_0$, then the integral in (\ref{eqn:g-int}) is $c_0$ for all
$z \in S_\eps$.  Hence, if $z_0 = x_0 + iy_0 \in S_\eps$ with $f(x) =
c_0$ for $x_0 - \delta \leq x \leq x_0 + \delta$, then using
(\ref{eqn:bound}), 
\begin{align*}
|g(z_0) - f(x_0)| &= \left| \frac{1}{\pi}\int_{-\infty}^\infty (f(x) -
c_0) \left(\frac{\eps}{(x-z)^2 + \eps^2}\right) dx \right| \\
&\leq \frac{C}{\pi}\left(\int_{-\infty}^{x_0-\delta} +
\int_{x_0+\delta}^\infty  \frac{\eps}{(x-x_0)^2} dx\right) \\
&\leq \frac{2\eps C}{\pi \delta}. 
\end{align*}

To show that ${\rm Re}\; g(z) \geq 0$, note that the second part of
(\ref{eqn:bound}) implies that ${\rm Re} \; (\eps/((x-z)^2 + \eps^2))
\geq 0$ for all $z \in S_\eps$, and since $f$ is real,
(\ref{eqn:g-int}) implies ${\rm Re} \; g(z) \geq 0$.  $\BOX$

\begin{lemma} \label{lemma:fb}
Let $V = \{(z,w): |w| < 1+|z|^2\}$.  Then there exists an injective
holomorphic map $\Phi: \CC^2 \ra V$.
\end{lemma}

{\bf Proof:} Let $H(z,w) = (w, w^2 - z/2)$.  Then $H$ is a polynomial
automorphism of $\CC^2$, and $(0,0)$ is an attracting fixed point for
$H$.  By \cite[appendix]{RR}, there is an injective holomorphic map
$\Psi$ from $\CC^2$ onto the basin of attraction of $(0,0)$, which is
defined as $B =\{p \in \CC^2: \lim_{n \ra \infty} H^n(p) = (0,0)\}$.
By \cite{BS}, there exists $R > 0$ such that $B$ is contained in 
$$ V_R = \{|z| \leq R, |w| < R\} \cup \{|z|\geq R, |w| < |z|\}. $$
Hence taking $\Phi = \Psi/R$ gives an injective holomorphic map from
$\CC^2$ into $V_1 \subseteq V$.  $\BOX$

\subsection{Proof of theorem \ref{thm:main}}

We will construct an automorphism of $\CC^2$ mapping $\Lambda_{0,r}$
into the complement of the set $V$ of lemma~\ref{lemma:fb}.  This will
be sufficient to prove the theorem, and by \cite{RR} this implies that
any discrete set contained in $\Lambda_{0,r}$ is tame.  

Choose an invertible, complex linear $A$ as in
lemma~\ref{lemma:discrete}.  Without loss of generality, we may
replace $\Lambda$ by $A(\Lambda)$, $p_j$ by $A(p_j)$, and $\Lambda_0$
by $A(\Lambda_0)$.  Then $\pi^1 \Lambda_0$ is contained in
$\cup_{k=1}^\infty L_k$, where each $L_k$ is a line of the form
$\RR+i\gamma_k$, $\gamma_k$ real.  Moreover, ${\rm dist}(L_j, L_k)
\geq 1$ if $j \neq k$, and $|p-q| \geq 1$ if $p,q \in \Lambda_0$
with $p \neq q$.  

Let $\{q_j\}_{j=1}^\infty$ be an enumeration of the set
$$ \{q \in \Lambda_0: |\pi^2 q| \leq 1/8\} = \{q \in \Lambda_0:
\ol{\Delta^2(q;1/8)} \cap (\CC \times \{0\}) \neq \emptyset\}. $$
Let $C = \log 32$, and define $f_k:\RR \ra [0,C]$ for each $k$ by
$$ f_k(x) = \left\{ \begin{array}{ll} 
0 & {\rm if} \ (x+i\gamma_k, 0) \in \ol{\Delta^2(q_j;1/8)} \ {\rm for
\ some} \ q_j \\ 
C & {\rm otherwise.}
\end{array} \right.
$$
Let $\delta = 1/16$, and choose $\eps \leq \delta/2$ small enough that
$2C\eps/\pi \delta \leq \log(3/2)$.  Let $r = \eps/2$, and recall that
$\Lambda_{0,r} = \cup_{p \in \Lambda_0} \Delta^2(p;r)$.  

\bigskip

Let $S_\eps^k = \{x+i(y+\gamma_k): -\eps < y < \eps\}$ and 
$U_\eps = \cup_{k=1}^\infty S_\eps^k$.
Define $g$ holomorphic on $U_\eps$ by applying lemma~\ref{lemma:g-int}
with $f=f_k$ to define $g$ on $S_\eps^k$.  By Arakelian's Theorem
(e.g. \cite{RR2}), there exists $h$ entire such that if $z \in
U_{\eps/2}$, then $|h(z) - g(z)| \leq \log(4/3)$.  Define
$$ F_1(z,w) = (z, w \exp(h(z))).  $$
Then $F_1: \CC^2 \ra \CC^2$ is biholomorphic.  

\bigskip

We show next that there is a complex line in the complement of
$F_1(\Lambda_{0,r})$.  To do this, let $p \in \Lambda_{0,r}$, and
suppose first that $p \in \Delta^2(q_j;r)$ for some $q_j$.  Choose $k$ so
that $\gamma_k = {\rm Im}\;\pi^1 q_j$, and write $\pi^1 p = x_0+iy_0$.  

Note that $|y_0 - \gamma_k| < r = \eps/2$.  Also, since $|\pi^2 q_j|
\leq 1/8$, we see that if $|x-x_0| < (1/8)-r$, then $(x+i\gamma_k, 0)
\in \ol{\Delta^2(q_j;1/8)}$.  Since $\delta<(1/8)-r$, we have
$f_k(x)=0$ for $x_0 - \delta \leq x \leq x_0 + \delta$, and hence by
lemma~\ref{lemma:g-int} and the choice of $\eps$ and $h$,
\begin{align*}
|h(\pi^1 p)| & \leq |g(\pi^1 p)| + \log (4/3) \\
& \leq \frac{2C \eps}{\pi \delta} + \log(4/3) \\
& \leq \log 2.
\end{align*}
Hence
\begin{equation} \label{eqn:in}
 |\pi^2 F_1(p)| \leq 2|\pi^2 p| \leq 2(|\pi^2 q_j| + r) \leq
\frac{1}{3}. 
\end{equation}

\bigskip

In the remaining case, $p \in \Lambda_{0,r}$ but $p \notin
\Delta^2(q_j;r)$ for any $j$, in which case $|\pi^2 p| \geq (1/8)-r$.
Let $q \in \Lambda_0$ such that $p \in \Delta^2(q;r)$, and choose $k$
so that $\gamma_k = {\rm Im}\;\pi^1 q$.  

Suppose first that $x_0 = {\rm Re} \; \pi^1 p$ satisfies $f_k(x) = C$ for
$|x-x_0| \leq \delta$.  Since $|y_0-\gamma_k| < r = \eps/2$, we have by
lemma~\ref{lemma:g-int} and choice of $\eps$ and $h$ that
\begin{align*}
{\rm Re}\; h(\pi^1 p) &\geq {\rm Re}\; g(\pi^1 p) - \log(4/3) \\
& \geq C - \frac{2C\eps}{\pi \delta} - \log(4/3) \\
& \geq \log 16. 
\end{align*}
Hence
\begin{equation} \label{eqn:out1}
 |\pi^2 F_1(p)| \geq 16|\pi^2 p| \geq 16((1/8)-r) >1. 
\end{equation}

Otherwise, $f_k(x) = 0$ for some $x$ with $|x-x_0| \leq \delta$, so
there exists $j$ such that $|\pi^1 p - \pi^1 q_j| \leq (1/8) +
\delta+r$, hence
$$ |\pi^1 q - \pi^1 q_j| \leq (1/8) + \delta + 2r \leq 1/4. $$
Since $q$ and $q_j$ are distinct points of $\Lambda_0$, we have
$|q-q_j| \geq 1$ by assumption, so $|\pi^2 q - \pi^2 q_j|^2 \geq
1-(1/4)^2$, and hence 
$$|\pi^2 q| \geq |\pi^2 q - \pi^2 q_j| - |\pi^2 q_j| \geq
\frac{\sqrt{15}}{4} - \frac{1}{8} $$
and 
$$ |\pi^2 p| \geq |\pi^2 q| - r \geq \frac{3}{4}. $$
Since ${\rm Re} \;g(\pi^1 p) \geq 0$ by lemma~\ref{lemma:g-int}, we have
${\rm Re}\;h(\pi^1 p) \geq -\log(4/3)$, and hence 
\begin{equation}  \label{eqn:out2}
 |\pi^2 F_1(p)| \geq \frac{3}{4} |\pi^2 p| \geq \frac{9}{16}. 
\end{equation}

{}From (\ref{eqn:in}), (\ref{eqn:out1}) and (\ref{eqn:out2}), we
conclude that if $p \in \Lambda_{0,r}$, then either $|\pi^2 F_1(p)|
\leq 1/3$ or $|\pi^2 F_1(p)| \geq 9/16$.  In particular,
$$ {\rm dist}(F_1(\Lambda_{0,r}), \CC \times \{\frac{1}{2}\})
\geq \frac{1}{16}.  $$
Note also that $\pi^1 F_1(p) = \pi^1 p$ for all $p \in \CC^2$.  

\bigskip

To finish the proof, we will construct $F_2$ similar to $F_1$ so that
$F_2(F_1(\Lambda_{0,r}))$ is contained in $\CC^2 \setminus V$, where
$V$ is as in lemma~\ref{lemma:fb}.

First note that for $z = x+iy \in S_\eps^k$, we have
$|y-i\gamma_k|<\eps$, so 
\begin{equation}  \label{eqn:re-g}
{\rm Re}\;[(z-i\gamma_k)^2 + (|\gamma_k|+r)^2 + 1 + \eps^2] \geq x^2 +
(|\gamma_k| + \eps)^2 + 1 > 0.
\end{equation}
Hence we can choose a branch of $\log$ so that 
\begin{equation}  \label{eqn:all-g}
g_2(z) = \log((z-i\gamma_k)^2 + (|\gamma_k|+r)^2 + 1 + \eps^2) +1 +
\log 16
\end{equation}
is holomorphic on $\cup_k S_\eps^k$.  Again by Arakelian's Theorem,
there exists $h_2$ entire such that if $z \in S_{\eps/2}^k$, then
$|g_2(z) - h_2(z)| \leq 1$, so ${\rm Re}\; h_2(z) \geq {\rm Re}\; g_2(z) -
1$.  Let
$$ F_2(z,w) = \left(z, \left(w-\frac{1}{2}\right) \exp(h_2(z))\right). $$
Again, $F_2:\CC^2 \ra \CC^2$ is biholomorphic.  Moreover, if $p \in
F_1(\Lambda_{0,r})$, then $|\pi^2 p - \frac{1}{2}| \geq 1/16$, and
$\pi^1 p = z = x+iy$ with $|y-\gamma_k| < r$ for some $k$, so by
(\ref{eqn:re-g}) and (\ref{eqn:all-g}), we have
\begin{align*}
|\pi^2 F_2(p)| & \geq \left|\pi^2p - \frac{1}{2}\right| 
\exp({\rm Re}\;h_2(z)) \\
& \geq \frac{1}{16} \exp({\rm Re}\;g_2(z) - 1) \\
& \geq x^2 + (|\gamma_k|+r)^2 + 1 \\
& \geq 1 + |\pi^1 p|^2 \\
& \geq 1 + |\pi^1 F_2(p)|^2.
\end{align*}

Hence $F_2 F_1(\Lambda_{0,r}) \cap V = \emptyset$, where $V \supseteq
\Phi(\CC^2)$ is as in lemma~\ref{lemma:fb}, so taking $F = F_1^{-1}
F_2^{-1} \Phi$ gives an injective holomorphic map $F:\CC^2 \ra
F_1^{-1} F_2^{-1}(V) \subseteq  \CC^2 \setminus \Lambda_{0,r}$ as
desired.  $\BOX$

\subsection{The general case of complements of small open balls}

It is now easy to deduce the following 
corollary from theorem \ref{thm:main}.

\begin{cor} Let $X$ be bimeromorphic to
a compact complex torus or to a Kummer K3 surface.
Then, given any finite set of points in $X$, the complement of a 
neighborhood of this set is dominable by $\CC^2$. In particular, 
such a complement is not measure hyperbolic.
\end{cor}

The case of elliptic fibrations over $\PP^1$ or over an elliptic curve
can be handled in the same way as that of theorem~\ref{thm:main}.
This is because removing
a finite number of small open balls (plus a smooth fiber away from them
if the base is $\PP^1$)
is tantamount to removing via the Jacobian fibration a discrete set
of contractible open sets in $\CC^2$ bounded away from the axis by
fixed constants and 
whose projection to the first factor $\CC$ is also a discrete
set of contractible open
subsets of $\CC$.  See also theorem~2.3 in \cite{Bu}.


\noindent
Gregery T. Buzzard\\
Department of Mathematics\\
Cornell University\\
Ithaca, NY  14853\\
USA

\medskip
\noindent
Steven Shin-Yi Lu\\
Department of Mathematics\\
University of Waterloo\\
Waterloo, Ontario, N2L3G1\\
Canada


\begin{thebibliography}{99}

 \bibitem[At]{At} M. F. Atiyah,
	{On analytic surfaces with double points,}
	{\it Proc. Royal. Soc. Lond., Ser A} 
	245 (1958), 237--244.

 \bibitem[Be]{Be}
	A. Beauville,
	{\it Complex algebraic surfaces},
	London Math. Soc. Lecture Note Series 68.

 \bibitem[Br]{Br}
	R. Brody,
	{Compact manifolds and hyperbolicity},
	{\it Trans. Amer. Math. Soc.},
	235 (1978), 213--219.

 \bibitem[Bu]{Bu}
	G. Buzzard,
	Tame sets, dominating maps, and complex tori,
	IHES preprint, IHES/M/99/46, 1999.


 \bibitem[BPV]{BPV}
	W. Barth, C. Peters, A. Van de Ven,
	{\it Compact Complex Surfaces},
	Springer-Verlag, Berlin, 1984.

 \bibitem[BM]{BM}
	E. Bierstone and P. Milman,
	Canonical desingularization in characteristic zero by blowing up 
	the maximal strata of a local invariant,
	{\it Invent. Math.}
	128 (1997), 207--302.

 \bibitem[BF]{BF}
	G. Buzzard and F. Forstneric,
	An interpolation theorem for holomorphic automorphisms of $\CC^n$,
	to appear in {\it J. Geom. Anal.}

 \bibitem[BS]{BS}
	E. Bedford and J. Smillie,
	Polynomial diffeomorphisms of $\CC^2$: currents,
	equilibrium measure and hyperbolicity,
	{\it Invent. Math.},
	103 (1991), no. 1, 69--99.

 \bibitem[CG]{CG}
	J. Carlson and P. Griffiths,
	A defect relation for equidimensional holomorphic mappings
        between algebraic varieties,
	{\it Ann. Math.}, (2) 95 (1972), 557--584.

 \bibitem[Del]{De} 
	P. Deligne, 
	Th\'eorie de Hodge I, 
	{\it Actes Cong. Int. Math.} (1970), 425--430.
	Th\'eorie de Hodge II,
	{\it Inst. Haut. Etud. Sci., Publ. Math.} 
	40 (1972), 5--57

 \bibitem[FK]{FK} 
	H. M. Farkas and I. Kra, 
	{\it Riemann Surfaces}, 
	Graduate Texts in Math. {71}, 
	Springer Verlag, New York, 1980.

 \bibitem[Gr]{green}
	M. Green,
	{\it Holomorphic maps to complex tori},
	Amer. J. Math., 100 (1978), no. 3, 615--620.
	
 \bibitem[GG]{GG}
	M. Green and P. Griffiths,
	Two applications of algebraic geometry to entire holomorphic mappings,
	{\it The Chern Symposium 1979}, 41--74 
	Springer-Verlag, New York Heidelberg-Berlin.

 \bibitem[GH]{GH}
	P. Griffiths and J. Harris,
	{\it Principles of algebraic geometry},	
	John Wiley and Sons, New York, 1978.

 \bibitem[Ha]{Ha}
	R. Hartshorne,
	{\it Algebraic Geometry},
	Graduate Texts in Math. 52,
	Springer Verlag, New York, 1977.

 \bibitem[Hi]{Hi}
	H. Hironaka,
	Resolution of singularities of an algebraic variety over a
	field of characteristic zero,
	{\it Ann. Math.},
	79 (1964) 109--326.

 \bibitem[Ii]{Ii}
	S. Iitaka,
	{\it Algebraic Geometry}, 
	Graduate Texts in Math. {76}, 
	Springer Verlag, New York, 1982.

 \bibitem[Ii1]{Ii1}
	S. Iitaka,
	Logarithmic forms of algebraic varieties,
	{\it J. Math. Soc. Japan}, 23 (1976), 525--544.

 \bibitem[In0]{In0}
	M. Inoue,
	On surfaces of class VII$_0$,
	{\it Invent. Math.} 24 (1974), 269--310.

 \bibitem[In]{In}
	M. Inoue,
	New surfaces with no meromorphic functions,
	{\it Proc. Int. Congr. Vancouver} 1974, 423--426,\\
	Ibid II, in {\it Complex Analysis and Algebraic Geometry},
	Iwanami-Shoten, Tokyo (1977), 91--106.

 \bibitem[K1]{K1}
 	Y. Kawamata,
	Addition formula of logarithmic Kodaira dimensions for morphisms of 
	relative dimension one, {\it Proceedings of the
	International Symposium on Algebraic Geometry at Kyoto in 1977},  
	Tokyo Kinokuniya (1978), 207--217.

 \bibitem[K2]{K2}
 	Y. Kawamata,
	On the classification of noncomplete algebraic surfaces, 
	Algebraic geometry (Proc. Summer Meeting, Univ.
	Copenhagen, Copenhagen, 1978), {\it Lecture Notes in Math.} (1979), 
	732,  215--232.

 \bibitem[K3]{K3}
	Y. Kawamata,
	Classification theory of non-complete algebraic surfaces, 
	{\it Proc. Japan Acad.} Ser. A Math. Sci. 54 (1978), no. 5,
	133--135.

 \bibitem[K4]{K4}
	Y. Kawamata,
	Characterization of abelian varieties,
	{\it Compositio Math.} 
	43 (1981), 253--276.

 \bibitem[Ko1]{K}
	K. Kodaira, Pluricanonical systems on algebraic surfaces of general
	type, {\it J. Math. Soc. Japan} 20 (1968), 170--192.

 \bibitem[Ko2]{Ko2}
	K. Kodaira, On compact complex analytic surfaces I,
	{\it Ann. Math.} 71 (1960), 111--152.
	II, {\it Ann. Math.} 77 (1963), 563--626.
	III,{\it Ann. Math.} 78 (1963), 1--40.

 \bibitem[Ko3]{Ko3}
	K. Kodaira, On the structure of compact complex analytic surfaces,
	Lecture Notes prepared in connection with the AMS Summer Institute
	on Algebraic Geometry held at the Whitney Estate, Woods Hole, Mass.
	July 1964. 

 \bibitem[Ko4]{Ko4}
	K. Kodaira, On the structure of compact complex analytic surfaces I,
	{\it Amer. J. Math.} 86 (1964), 751--798.

 \bibitem[KO]{KO}
	S. Kobayashi and T. Ochiai,
	Meromorphic mappings onto compact complex spaces of general type,
	{\it Invent. Math.},
	31, (1975), 7--16.

 \bibitem[Kob]{Kob}
	S.  Kobayashi,
	{\it Hyperbolic manifold and holomorphic mappings},
	Marcel Decker, New York, 1970.


 \bibitem[Lang]{Lang}
	S. Lang,
	Hyperbolic and Diophantine analysis,
	{\it Bull. Amer. Math. Soc.} 
	14,(1986) no. 2, 159--205.

 \bibitem[LP]{LP}
	E. Looijenga, C. Peters,
	Torelli theorems for K\"ahler K3 surfaces,
	{\it Compositio Mathematica},
	42 (1981), 145--186.

 \bibitem[Lu1]{Lu1}
	S. Lu,
	On meromorphic maps into varieties of log-general type,
	{\it Proceedings of Symposia in Pure Mathematics},
	52 (1991), part 2, 305--333.

 \bibitem[M]{M}
	M. Miyanishi,  
	{\it Non-complete algebraic surfaces},
	Lecture Notes in Math. 857,
	Springer-Verlag 1981.

 \bibitem[Mas]{Mas}
	B. Maskit,
	{\it Kleinian Groups}, G.M.W 287
	Springer-Verlag 1988.

 \bibitem[MM]{MM} S. Mori, S. Mukai,
	The uniruledness of the moduli space of curves of genus 11,
	{\it Lecture Notes in Math.},
	1019 (1982), 334--353.


 \bibitem[PS]{PS}	
	I.I. Piateckii-Shapiro, I.R. Shafarevic,
	A Torelli theorem for algebraic surfaces of type K-3,
	{\it Izv. Akad. Nauk. SSSR},
	Ser. Math. 35 (1971), 503--572.
	

 \bibitem[RR1]{RR}
	J.-P. Rosay and W. Rudin,
	Holomorphic maps from $\CC^n$ to $\CC^n$,
	{\it Trans. Amer. Math. Soc.},
	310 (1988), no. 1, 47--86.

 \bibitem[RR2]{RR2}
	J.-P. Rosay and W. Rudin,
	Arakelian's approximation theorem, 
	{\it Amer. Math. Monthly}, 
	96 (1989), no. 5, 432--434.

 \bibitem[R]{rudin}
	W. Rudin,
	{\it Real and complex analysis},
	McGraw Hill, New York, 1987.

 \bibitem[Shi]{Shi} T. Shioda, 
	The period map of abelian surfaces,
	{\it J. Fac. Sci. Univ. Tokyo,} Sect. IA,
	25 (1978), 47--59.

 \bibitem[Siu]{Siu} Y. T. Siu,
	Every K3-surface is K\"ahler,
	{\it Invent. Math.},
	73 (1983), 139--150.

 \bibitem[T]{T}
	A. Todorov,
	Applications of the K\"{a}hler-Einstein-Calabi-Yau metric to
	moduli of K3-surfaces,
	{\it Invent. Math.},
	8 (1980), 251--255.

 \bibitem[Ue]{Ue}
	K. Ueno,
	{\it Classification theory of algebraic varieties and compact complex
	Spaces}, 1975, Springer-Verlag, New York.

 \bibitem[Yau]{Yau}
	S. T. Yau,
	On the Ricci-curvature of a complex K\"ahler manifold and
	the complex Monge-Amp\`ere equation,
	{\it Comm. Pure Appl. Math.},
	31 (1978), 339--411.

\end{thebibliography}
\end{document}